\documentclass[a4paper,10pt]{article}

\usepackage{lineno,hyperref}
\modulolinenumbers[5]

\usepackage{enumitem}
\usepackage{amsmath,dsfont,amsthm}
\usepackage{amssymb}
\usepackage{adjustbox,lipsum}
\usepackage{mathtools}
\usepackage{fancyhdr}
\usepackage{bibentry}
\usepackage{graphicx}
\usepackage{xcolor}
\usepackage{stackengine}
\allowdisplaybreaks

\newtheorem{theorem}{Theorem}[section]
\newtheorem{lemma}[theorem]{Lemma}
\newtheorem{definition}[theorem]{Definiton}
\newtheorem{proposition}[theorem]{Proposition}

\newtheorem{assumption}[theorem]{Assumption}

\theoremstyle{definition}
\newtheorem{remx}[theorem]{Remark}
\newenvironment{remark}
  {\pushQED{\qed}\remx}
  {\popQED\endremx}

\makeatletter
\newsavebox\myboxA
\newsavebox\myboxB
\newlength\mylenA

\newcommand*\yoverline[2][0.75]{%
    \sbox{\myboxA}{$\m@th#2$}%
    \setbox\myboxB\null
    \ht\myboxB=\ht\myboxA%
    \dp\myboxB=\dp\myboxA%
    \wd\myboxB=#1\wd\myboxA
    \sbox\myboxB{$\m@th\overline{\copy\myboxB}$}
    \setlength\mylenA{\the\wd\myboxA}
    \addtolength\mylenA{-\the\wd\myboxB}%
    \ifdim\wd\myboxB<\wd\myboxA%
       \rlap{\hskip 0.5\mylenA\usebox\myboxB}{\usebox\myboxA}%
    \else
        \hskip -0.5\mylenA\rlap{\usebox\myboxA}{\hskip 0.5\mylenA\usebox\myboxB}%
    \fi}
\makeatother
\numberwithin{equation}{section}
\begin{document}




\newcommand{\diver}{\operatorname{div}}
\newcommand{\ww}{{W}}
\newcommand{\yy}{{y}}
\newcommand{\lin}{\operatorname{Lin}}
\newcommand{\curl}{\operatorname{curl}}
\newcommand{\ran}{\operatorname{Ran}}
\newcommand{\kernel}{\operatorname{Ker}}
\newcommand{\la}{\langle}
\newcommand{\ra}{\rangle}
\newcommand{\N}{\mathbb{N}}
\newcommand{\R}{\mathbb{R}}
\newcommand{\ld}{\lambda}
\newcommand{\fai}{\varphi}
\newcommand{\hh}{H^1(M)}
\newcommand{\hhh}{H^1_0(M)}
\newcommand{\bsigma}{\boldsymbol{\sigma}}
\newcommand{\D}{\boldsymbol{\mathrm{D}}}
\newcommand{\C}{{\mathds{C}}}
\newcommand{\e}{\mathbf{e}}
\newcommand{\E}{{E}}
\newcommand{\eps}{{\boldsymbol{\epsilon}}}
\newcommand{\vareps}{\boldsymbol{{\varepsilon}}}
\newcommand{\Bsigma}{\boldsymbol{{\Sigma}}}
\newcommand{\PP}{{\boldsymbol{\mathrm{P}}}}
\newcommand{\lol}{{\ell}}
\newcommand{\BPP}{\bar{{\mathbf{P}}}}
\newcommand{\TT}{\boldsymbol{t}}
\newcommand{\ta}{\mathrm{t}_{\bsigma}}
\newcommand{\tb}{\mathrm{t}_{\D}}
\newcommand{\tc}{\mathrm{t}_{\PP}}
\newcommand{\fa}{\mathrm{f}_{\bsigma}}
\newcommand{\fb}{\mathrm{f}_{\D}}
\newcommand{\fs}{\boldsymbol{f_{\boldsymbol{\sigma}}}}
\newcommand{\fd}{\boldsymbol{f_{\D}}}
\newcommand{\fp}{\boldsymbol{f_{\PP}}}
\newcommand{\fc}{\mathrm{f}_{\PP}}
\newcommand{\pai}{\boldsymbol{\pi}}
\newcommand{\0}{0}
\newcommand{\n}{\mathbf{n}}
\newcommand{\uu}{{\boldsymbol{\mathrm{u}}}}
\newcommand{\buu}{\bar{{\mathbf{u}}}}
\newcommand{\bPP}{\bar{{P}}}
\newcommand{\bphi}{\bar{\phi}}
\newcommand{\x}{{x}}
\newcommand{\y}{{y}}
\newcommand{\sss}{{S}}
\newcommand{\Omaga}{{\Omega}}
\newcommand{\bs}{\boldsymbol}
\newcommand{\q}{{\boldsymbol{\mathrm{q}}}}
\newcommand{\ii}{\mathcal{I}}
\newcommand{\ee}{\mathcal{E}}
\newcommand{\bolde}{{E}}
\newcommand{\gggg}{{\Gamma}}
\newcommand{\vv}{{\mathbf{v}}}
\newcommand\barbelow[1]{\stackunder[1.2pt]{$#1$}{\rule{.8ex}{.075ex}}}
\newcommand{\pk}{\boldsymbol{P}_{k}^{\tau}}
\newcommand{\pko}{\boldsymbol{P}_{k+1}^{\tau}}
\newcommand{\tk}{t_{k}^{\tau}}
\newcommand{\tko}{t_{k+1}^{\tau}}
\newcommand{\tauk}{\tau_k}
\newcommand{\tauko}{\tau_{k+1}}
\newcommand{\ptauu}{\bar{\PP}_{\tau}}
\newcommand{\ptaud}{\barbelow{\PP}_{\tau}}
\newcommand{\ptauh}{\hat{\PP}_{\tau}}
\newcommand{\ttauu}{\bar{t}_{\tau}}
\newcommand{\ttaud}{\barbelow{t}_{\tau}}
\newcommand{\ptauuj}{\bar{\PP}_{\tau^j}}
\newcommand{\ptaudj}{\barbelow{\PP}_{\tau^j}}
\newcommand{\ptauhj}{\hat{\PP}_{\tau^j}}
\newcommand{\ttauuj}{\bar{t}_{\tau^j}}
\newcommand{\ttaudj}{\barbelow{t}_{\tau^j}}
\newcommand{\nnn}{\boldsymbol{n}}
\newcommand{\fff}{{\ell}}
\newcommand{\gog}{\boldsymbol{g}}
\newcommand{\puu}{\Gamma_{\uu}}
\newcommand{\pphi}{\Gamma_{\phi}}
\newcommand{\pPP}{\Gamma_{\PP}}
\newcommand{\FFF}{{f}}
\newcommand{\gtwos}{\mathcal{G}_{2,grad,s}}
\newcommand{\sos}{{F}}
\newcommand{\LL}{\mathrm{L}}
\newcommand{\W}{\mathbb{W}_p}
\newcommand{\WW}{\mathbb{W}_{p'}}
\newcommand{\V}{\mathbb{V}_{p,r}}
\newcommand{\U}{\mathbb{U}}
\newcommand{\SSS}{\mathcal{S}}
\newcommand{\ten}{\\[1.8pt]}
\newcommand{\six}{\\[1.8pt]}
\newcommand{\nb}{\nonumber}

\title{On the local in time well-posedness of an elliptic-parabolic ferroelectric phase-field model}
\author{Yongming Luo \thanks{Institut f\"{u}r Wissenschaftliches Rechnen, Technische Universit\"at Dresden, Germany} \thanks{\href{mailto:yongming.luo@tu-dresden.de}{Email: yongming.luo@tu-dresden.de}}
}

\date{}
\maketitle

\begin{abstract}
We consider a state-of-the-art ferroelectric phase-field model arising from the engineering area in recent years, which is mathematically formulated as a coupled elliptic-parabolic differential system. We utilize the maximal parabolic regularity theory to show the local in time well-posedness of the ferroelectric problem in both 2D and 3D spaces, which is sharp in the sense that the local solution is unique and a blow-up criterion is present. The well-posedness result will firstly be proved under some general assumptions. Afterwards we give sufficient geometric and regularity conditions which will guarantee the fulfillment of the imposed assumptions.
\end{abstract}

\section{Introduction}
In this paper, we study a state-of-the-art ferroelectric phase-field model \cite{Schrade2013,GAMM2015} arising from the engineering area in recent years. To formulate the problem, we list the corresponding physical quantities in Table \ref{table of var} below.

Given a bounded domain $\Omaga\subset\R^d$ and a time slot $(0,T)$ with $T\in(0,\infty)$, the following constitutive laws are imposed for $\bsigma$ and $\D$:
\begin{alignat}{2}
\bsigma&=\C(\PP)\big(\vareps(\uu)-\vareps^0(\PP)\big)+\e(\PP)^T\nabla\phi &&\quad\quad\text{ in }(0,T)\times\Omaga,\\%
\D&=\e(\PP)\big(\vareps(\uu)-\vareps^0(\PP)\big)-\eps(\PP)\nabla\phi+\PP &&\quad\quad\text{ in }(0,T)\times\Omaga.
\end{alignat}
\begin{table}
\caption{List of variables}
\label{table of var}       
\renewcommand{\arraystretch}{1.2}
\begin{tabular}{|r|c|}
\noalign{\vspace*{8pt}}\hline
Mechanical displacement& $\uu:[0,T]\times\Omaga\to\R^d$\\
Infinitesimal small strain tensor &$\vareps(\uu):= (\nabla \uu+\nabla \uu^T)/2$\\
Electric potential&$\phi:[0,T]\times\Omaga\to\R$\\
Spontaneous polarization &$\PP:[0,T]\times\Omaga \to  \R^d$\\
Cauchy stress tensor &$\bsigma:[0,T]\times\Omaga \to\mathrm{Lin}_{\mathrm{sym}}(\R^d)$\\
Dielectric displacement &$\D:[0,T]\times\Omaga \to\R^d$\\
Elastic stiffness tensor &${\C}:\R^d \to\mathrm{Lin}_{\mathrm{sym}}(\R^{d\times d})$ with $\C_{ij}^{kl}=\C_{ji}^{lk}$\\
Symmetric plastic strain tensor &$\vareps^0:\R^d  \to \mathrm{Lin}_{\mathrm{sym}}(\R^d)$\\
Coupling effect tensor &$\e:\R^d  \to\mathrm{Lin}(\R^{d\times d},\R^d)$\\
Symmetric dielectric matrix &$\eps:\R^d  \to\mathrm{Lin}_{\mathrm{sym}}(\R^d)$\\
Separation energy density &$\omega:\R^d  \to\R$\\
\hline\noalign{\smallskip}
\end{tabular}
\end{table}

\noindent We aim to find a solution $(\uu,\phi,\PP)$ of the Cauchy problem
\begin{subequations}\label{weakdepa}
\begin{alignat}{2}
-\diver \bsigma=\fa& &&\quad\quad\text{ in }(0,T)\times\Omaga,\label{cauchy momentum law}\\
-\diver \D=\fb& &&\quad\quad\text{ in }(0,T)\times\Omaga,\label{gauss law}\\
\PP'=\Delta\PP-D_{\PP}H(\vareps(\uu),\nabla\phi,\PP)-D_{\PP}\omega(\PP)+\fc& &&\quad\quad\text{ in }(0,T)\times\Omaga\label{2nd law of thermodynamics}
\end{alignat}
\end{subequations}
with the boundary conditions
\begin{subequations}\label{weakdepb}
\begin{alignat}{2}
\uu&=\0&&\quad\text{ on }(0,T)\times\Gamma_{\uu},\\
\phi&=0&&\quad\text{ on }(0,T)\times\Gamma_{\phi},\\
\PP&=\0&&\quad\text{ on }(0,T)\times\Gamma_{\PP},\label{ppd}\\
\bsigma \n&=\ta&&\quad\text{ on }(0,T)\times\Gamma^\mathrm{N}_{\uu},\\
\D\cdot \n&=\tb&&\quad\text{ on }(0,T)\times\Gamma^\mathrm{N}_{\phi},\\
\partial\PP/\partial\n&=\tc&&\quad\text{ on }(0,T)\times\Gamma^\mathrm{N}_{\PP}\label{ppn}
\end{alignat}
\end{subequations}
and the initial condition
\begin{alignat}{2}\label{initial}
\PP(0,{\x})=\PP_0({\x})\quad\text{ in }\Omaga,
\end{alignat}
where $\Gamma_{\uu}\,\dot{\cup}\,\Gamma^\mathrm{N}_{\uu}=\Gamma_{\phi}\,\dot{\cup}\,\Gamma^\mathrm{N}_{\phi}=\Gamma_{\PP}\,\dot{\cup}\,\Gamma^\mathrm{N}_{\PP}=\partial\Omaga$\footnote{Heuristically, $\Gamma_{\uu},\Gamma_{\phi},\Gamma_{\PP}$ and $\Gamma^\mathrm{N}_{\uu},\Gamma^\mathrm{N}_{\phi},\Gamma^\mathrm{N}_{\PP}$ denote the Dirichlet and Neumann boundary portions respectively.}, $\n$ is the outer normal vector, $\fa,\fb,\fc$ are the body force densities, $\ta,\tb,\tc$ are the boundary flux densities and $H$ is the bulk energy density defined by
\begin{align}\label{H}
H(\vareps,\nabla\phi,\PP)=\frac{1}{2}\mathbb{B}(\PP)\begin{pmatrix}\vareps-\vareps^0(\PP)\\ -\nabla\phi\end{pmatrix}:\begin{pmatrix}\vareps-\vareps^0(\PP)\\ \nabla\phi\end{pmatrix}+\nabla\phi\cdot\PP
\end{align}
with
\begin{align}\label{not symmetric tentor}
\mathbb{B}(\PP)=\begin{pmatrix}\C(\PP)&\e(\PP)^T\\ -\e(\PP)&\eps(\PP)\end{pmatrix}.
\end{align}
A schematic diagram for the boundary partition $\partial\Omaga=\Gamma_{\uu}\,\dot{\cup}\,\Gamma^\mathrm{N}_{\uu}$ is depicted in Fig. \ref{boundary} below.
\begin{figure}[!htbp]
  \centering
  \includegraphics[width=40mm]{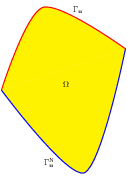}
  \caption{An illustration for the boundary partition $\partial\Omaga=\Gamma_{\uu}\,\dot{\cup}\,\Gamma^\mathrm{N}_{\uu}$. The red part and blue part are the Dirichlet and Neumann boundary portions corresponding to the variable $\uu$ respectively.}\label{boundary}
\end{figure}

The study of ferroelectric materials has a longstanding history since its very first discovery by Valasek \cite{SchroderRoman2005} in the 20th century. A crucial property of the ferroelectric materials is their ability of transforming an electric signal into mechanical deformation, and vice versa. Such property, which is known as the piezoelectricity, displays the linear counterpart of a ferroelectric model. However, when the ferroelectric materials are cooled down below the Curie temperature, they enter a nonlinear regime and their dynamics are described by a nonlinear evolutionary law. Such phase transition phenomenon also changes the dipole moments between the electric charges within the single crystals, hence also leads to the spontaneous polarization. Moreover, the spontaneous polarization is aligned along or opposite to the axes of the single tetrahedral crystals, therefore $90^{\circ}$ or $180^{\circ}$ domain walls are also obtained. Each domain wall can thus be seen as a phase-field parameter, and this important feature enables us to study the ferroelectric model using the phase-field theory. The above mentioned properties particularly reveal the suitability of ferroelectric materials for applications involving memory effects, such as sensors, actuators, capacitors etc. For a more comprehensive introduction on ferroelectric materials, we refer to \cite{Werner1964,Cross1987,Haertling1999,Huber1999,Kamlah2001} and the references therein. We also refer the readers to \cite{Bartic2001,SchroderRoman2005,SuLandis2007,Marton2010,Kaltenbacher2011,Miehe2011,GAMM2015} for recent progress on the modeling and numerical simulations of the ferroelectric model.

In contrary to the diverse results for ferroelectrics given in the engineering field, a rigorous mathematical qualitative or quantitative analysis for the ferroelectric phase-field model can barely be found in the literature. The first remarkable result concerning the well-posedness of the ferroelectric phase-field model was given by Mielke and Timofte \cite{MielkeTimofte2006}. Therein, the ferroelectric problem is considered on a bounded domain $\Omaga\subset\R^d$ for a time slot $[0,T]$ and formulated in terms of an abstract Helmholtz free energy $\mathcal{H}$ and a non-negative convex dissipation functional $\Psi$. The key assumption made in \cite{MielkeTimofte2006} is the \textbf{rate-independence} of the dissipation functional $\Psi$, i.e.
$$ \forall\,r\geq 0\,\,\forall\,\PP\in{\boldsymbol{\mathrm{X}}}_3:\,\Psi(r\PP)=r\Psi(\PP),$$
where ${\boldsymbol{\mathrm{X}}}_3$ is some suitable underlying function space for $\PP$. Relying on certain generalized De Giorgi energy-dissipation principles for rate-independent systems (see for instance \cite{MielkeRoubicek2015}), Mielke and Timofte were able to show the existence of a so called energetic solution for rate-independent ferroelectric system. Similar models are considered in \cite{KutterSaendig2011} and \cite{Pechstein2020}, where the conditions on $\mathcal{H}$ are different to the ones in \cite{MielkeTimofte2006} (more precisely, there is no gradient energy term $\frac{1}{2}\|\PP\|^2_{L^2(\Omaga)}$ contained in $\mathcal{H}$, and the evolutionary law for $\PP$ reduces to an ODE) but the dissipation functionals $\Psi$ remain rate-independent. Recently, different solution schemes such as BV-solutions and vanishing viscosity solutions etc. have been successively applied for rate-independent ferroelectric phase-field models. We refer to \cite{Knees2019,LuoThesis2019} for further details.

In general, however, ferroelectric materials display \textbf{rate-dependence} on the spontaneous polarization, and the system can be assumed rate-independent only when the external loadings possess a relatively much lower frequency compared to the internal dipoles. In this paper, we will study the rate-dependent ferroelectric model introduced in \cite{Schrade2013,GAMM2015}, which is considered as one of the state-of-the-art models in the current research for ferroelectric materials. We focus on the 2D and 3D cases, being the general scenarios under consideration in practical applications. To the author's best knowledge, this is the first rigorous well-posedness result for rate-dependent ferroelectric phase-field models with purely quadratic dissipation growth. A number of interesting features are also displayed by the model:\begin{itemize}
\item The spontaneous polarization $\PP$ admits boundary conditions of mixed type\footnote{More precisely, $\PP$ satisfies  \eqref{ppd} and \eqref{ppn} with $\partial\Omega=\Gamma_{\PP}\,\dot{\cup}\,\Gamma^\mathrm{N}_{\PP}$.}, where $\PP$ in the models \cite{MielkeTimofte2006,Knees2019} with pure Neumann boundary conditions are considered.
\item The material parameters ${\C},\vareps^0, \e,\eps$ are assumed to be functions of $\PP$, which are usually in literature merely assumed to be (fully or partially) constants \cite{KutterSaendig2011,Schrade2013} or depend on the spatial variable ${\x}\in\Omaga$ \cite{Knees2019}.
\item When the variable $(\uu,\phi)$ possesses pure Dirichlet boundary conditions, then no uniform ellipticity and boundedness conditions are imposed on the material parameters, see Theorem \ref{cube or C^1 domain} below.
\item The well-posedness result is sharp, in the sense that the local solution is unique and a blow-up criterion is present, see Definition \ref{solution concept} below.
\end{itemize}
We briefly summarize the idea for proving the local well-posedness of the model problem. The proof follows the same lines of \cite{MeinschmidtMeyerRehberg16}: by multiplying a test function and using partial integration, the balance law \eqref{weakdepa} is reformulated to the following coupled elliptic-parabolic system:
\begin{subequations}\label{elliptic-parabolic-1}
\setlength{\jot}{10pt}
\begin{alignat}{2}
-\mathrm{Div}\left(\mathbb{B}(\PP(t))\begin{pmatrix}\vareps\big(\uu(t)\big) \\ \nabla\phi(t)\end{pmatrix}\right)=\lol(t,\PP(t))& &&\quad\text{ in }\W^*,\label{elliptic-parabolic-1a} \\
\PP'(t)=\Delta\PP(t)+\sos(t,\vareps(\uu(t)),\nabla\phi(t),\PP(t))& &&\quad\text{ in }\big(W^{-1,p}_{\Gamma_{\PP}}(\Omaga)\big)^d,
\end{alignat}
\end{subequations}
where the space $\W^*$ is defined by \eqref{special not2}. Using Lax-Milgram we infer that $(\uu(t),\phi(t))$ as well as the nonlinear functional $\sos(t,\vareps(\uu(t)),\nabla\phi(t),\PP(t))$ are uniquely determined by the given pair $(t,\PP(t))$. Hence $\sos(t,\vareps(\uu(t)),\nabla\phi(t),\PP(t))$ can be rewritten to the form $S(t,\PP(t))$ and \eqref{elliptic-parabolic-1} equivalently reduces to a semilinear parabolic equation
\begin{align}\label{single parabolic-1}
\PP'(t)=\Delta\PP(t)+S(t,\PP(t)).
\end{align}
We then exploit the existence and uniqueness result given by Cl\'{e}ment and Li \cite{LiClement1994} to obtain the local well-posedness of the semilinear parabolic problem \eqref{single parabolic-1}. It is worth noting that the key point for the proof is that the solution $(\uu,\phi)$ of \eqref{elliptic-parabolic-1a} is of class $W^{1,p}$ for some $p>d$ (see \eqref{why p larger d} below). This property will firstly be imposed as an assumption (Assumption \ref{A3}) for the local well-posedness result. We then give sufficient geometric and regularity conditions which will guarantee the validity of Assumption \ref{A3}. For details, we refer to Section \ref{main results}.

\subsection{Outline of the paper}
The rest of the paper is organized as follows: In Section \ref{main results} we state the precise statements of the main results; In Section \ref{notations} we list the notation and definitions which will be used throughout this paper; finally, Section \ref{proof of thm 1} is devoted to the proof of Theorem \ref{C^1 domain} and Section \ref{proof of thm 2} to the proofs of Theorems \ref{2d mixed}, \ref{cube or C^1 domain} and \ref{poly thm}.
\section{Assumptions and main results}\label{main results}
In this section we state the main results of the paper. For the corresponding notation and definitions, we refer to Section \ref{notations}. We firstly define the precise solution concept for the Cauchy problem \eqref{weakdepa} to \eqref{initial}:
\begin{definition}\label{solution concept}
Let $p\in(d,\infty)$ and $r\in (\frac{2p}{p-d},\infty)$. We say that $(\uu,\phi,\PP)$ with $\PP(0,{\x})=\PP_0({\x})$ is a solution of the Cauchy problem \eqref{weakdepa} to \eqref{initial} on the maximal time interval $(0,\hat{T})$ with $\hat{T}\in(0,T]$ if for each $S\in(0,\hat{T})$, $(\uu,\phi,\PP)$ satisfies the equations
\begin{subequations}\setlength{\jot}{10pt}
\begin{alignat}{2}
-\mathrm{Div}\left(\mathbb{B}(\PP(t))\begin{pmatrix}\vareps\big(\uu(t)\big) \\ \nabla\phi(t)\end{pmatrix}\right)=\lol(t,\PP(t))& &&\quad\text{ in }\W^*,\label{elliptic part}\\
\PP'(t)=\Delta\PP(t)+\sos(t,\vareps(\uu(t)),\nabla\phi(t),\PP(t))& &&\quad\text{ in }\big(W^{-1,p}_{\Gamma_{\PP}}(\Omaga)\big)^d\label{parabolic part}
\end{alignat}
\end{subequations}
for a.a. $t\in(0,S)$ and possesses the regularity
\begin{gather}
(\uu,\phi)\in L^{2r}(0,S;\W),\\[4pt]
\PP\in W^{1,r}(0,S;\big(W^{-1,p}_{\Gamma_{\PP}}(\Omaga)\big)^d)\cap L^{r}(0,S;\big(W^{1,p}_{\Gamma_{\PP}}(\Omaga)\big)^d).\label{solution3d}
\end{gather}
Here, the spaces $\W$ and $\W^*$ are defined by
\begin{align}
\W&:=\big(W^{1,p}_{\Gamma_{\uu}}(\Omaga)\big)^d\times W^{1,p}_{\Gamma_{\phi}}(\Omaga),\\
\W^{*}&:=\big(W^{-1,p}_{\Gamma_{\uu}}(\Omaga)\big)^d\times W^{-1,p}_{\Gamma_{\phi}}(\Omaga).
\end{align}
Moreover, the maximality of $\hat{T}$ is understood as follows:
\begin{itemize}
\item[(a)] Either \eqref{solution3d} does not hold for $S=\hat{T}$, or
\item[(b)] $\hat{T}=T$ and \eqref{solution3d} also holds for $S=T$.
\end{itemize}
In the cases (a) and (b), the solution is said to be a finite time blow-up and a global solution respectively. Finally, the tensor $\mathbb{B}(\PP)$ is defined by \eqref{not symmetric tentor} and the functionals $\lol(t,\PP),\,\sos(t,\vareps,\nabla\phi,\PP)$ are defined by
\begin{align}
\lol(t,\PP):&=\lol_1(\PP)+\lol_2(t),\label{ltp}\\
\nonumber\\
\sos(t,\vareps,\nabla\phi,\PP):&=-D_{\PP}H(\vareps,\nabla\phi,\PP)-D_{\PP}\omega(\PP)+\lol_3(t)\label{stp intro},\\
\nonumber\\
\lol_1(\PP)[\uu,\phi]:&=\begin{pmatrix}-\mathrm{Div}_{\uu}\left(\C(\PP)\vareps^0(\PP)\right)\\ \mathrm{Div}_{\phi}\left(\e(\PP)\vareps^0(\PP)-\PP\right)\end{pmatrix}\left[\uu,\phi\right]\nonumber\\
&=\int_{\Omaga}\C(\PP)\vareps^0(\PP):\vareps(\uu)-(\e(\PP)\vareps^0(\PP)-\PP)\cdot\nabla\phi\,d\x,\label{ltp1}
\end{align}
\begin{align}
\lol_{2}(t)[\uu,\phi]:&=\int_{\Omaga}\fa(t,\x)\cdot \uu(\x)-\fb(t,\x) \phi(\x)\,d\x+\int_{\Gamma^\mathrm{N}_{\uu}}\ta(t,\x)\cdot\uu(\x)\,d\sss\nonumber\\
&-\int_{\Gamma^\mathrm{N}_{\phi}}\tb(t,\x)\phi(\x)\,d\sss,\label{ltp2}\\
\nonumber\\
\lol_3(t)[\PP]:&=\int_{\Omaga}\fc(t,\x)\cdot \PP(\x)\,d\x+\int_{\Gamma^\mathrm{N}_{\PP}}\tc(t,\x)\cdot\PP(\x)\,d\sss
\end{align}
for a test function $(\uu,\phi,\PP)\in\mathbb{W}_{p'}\times \big(W^{1,p'}_{\Gamma_{\PP}}(\Omaga)\big)^d$.
\end{definition}
\subsection{Assumptions}
We impose the following assumptions for the well-posedness result. In particular, the somehow indirectly formulated Assumption \ref{A3} plays a key role for establishing the maximal parabolic regularity of the Laplacian $-\Delta: W_{\Gamma_{\PP}}^{1,p}(\Omaga)\to W_{\Gamma_{\PP}}^{-1,p}(\Omaga)$. In general, it is not an easy task to verify the validity of the abstract Assumption \ref{A3}. We will give in Section \ref{proof of thm 2} precise conditions, under which the Assumption \ref{A3} is satisfied.
\begin{enumerate}[label=(A\arabic*)]
\item \label{A1} $\Omega\subset\R^d$, $d\in\{2,3\}$, is a bounded domain, $\Gamma_{\PP}$ is $(d-1)$-set with non-zero surface measure and $\Omega\cup\Gamma_{\PP}$ is Gr\"oger-regular (see Definition \ref{def of geometric sets}).

\item \label{A2}$\C,\e,\vareps^0,\eps,\omega$ are $C^{1,1}$-functions on $\R^d$.

\item \label{A3}
There exists some $p\in I_d$, where
\begin{align}
I_d=\left\{
             \begin{array}{ll}
             (2,\infty), &\text{if\ }d=2,  \\
             (3,6], &\text{if\ }d=3,
             \end{array}
\right.
\end{align}
such that the linear differential operator $\LL_{\PP}$ given by
\begin{align*}
\LL_{\PP}(\uu,\phi):=-\mathrm{Div}\left(\mathbb{B}(\PP)\begin{pmatrix}\vareps\big(\uu\big) \\ \nabla\phi\end{pmatrix}\right)
\end{align*}
defines a topological isomorphism from $\W$ to $\W^*$ for all $\PP$ in the space $\U$ defined by \eqref{special not4}.

\item \label{A4} There exists some $r\in(\frac{2p}{p-d},\infty)$, where $p$ is the number given by Assumption \ref{A3}, such that
\begin{align*}
\lol_2&\in L^{2r}(0,T;\W^* ),\nonumber\\
\lol_3&\in L^{r}(0,T;\big(W_{\Gamma_{\PP}}^{-1,p}(\Omaga)\big)^d).\nonumber
\end{align*}

\end{enumerate}


\begin{remark}
The restriction $p\in I_d$ given in Assumption \ref{A3} comes from the fact that in general it is only expected that the Laplacian $-\Delta: W_{\Gamma_{\PP}}^{1,p}(\Omaga)\to W_{\Gamma_{\PP}}^{-1,p}(\Omaga)$ admits the maximal parabolic regularity for $p\in I_d$, see \cite[Thm. 1.6]{Egert2018}.
\end{remark}
\subsection{Main results}
Having given the Assumptions \ref{A1} to \ref{A4} we are ready to formulate the local in time well-posedness result:
\begin{theorem}\label{C^1 domain}
Let the Assumptions \ref{A1} to \ref{A4} be satisfied. Assume also that
$$\PP_0\in(\big(W^{1,p}_{\Gamma_{\PP}}(\Omaga)\big)^d,\big(W^{-1,p}_{\Gamma_{\PP}}(\Omaga)\big)^d)_{\frac{1}{r},r},$$
where $(\cdot,\cdot)_{\frac{1}{r},r}$ denotes the real interpolation. Then the Cauchy problem \eqref{weakdepa} to \eqref{initial} has a unique local solution $(\uu,\phi,\PP)$ in the maximal time interval $(0,\hat{T})$ for some $0<\hat{T}\leq T$ in the sense of Definition \ref{solution concept}.
\end{theorem}
In practical applications, the fulfilment of the Assumptions \ref{A1}, \ref{A2} and \ref{A4} can be achieved by imposing sufficiently smooth conditions on the underlying domain and the external loadings. The problem hence reduces to finding precise conditions, under which the abstractly formulated Assumption \ref{A3} is satisfied. For 2D case, the Assumption \ref{A3} is an immediate consequence of the elliptic regularity result given by \cite{HDJonssonKneesRehberg2016} (Theorem \ref{higherorder} below). More precisely, we have the following:
\begin{theorem}\label{2d mixed}
Suppose that the following assumptions are satisfied:
\begin{enumerate}[label=(B\arabic*)]
\item \label{B1} $\Omaga\subset\R^2$ is a bounded domain, $\Gamma_{\uu}$, $\Gamma_{\phi}$ are closed $1$-sets with non-zero surface measure and $\Omega\cup \Gamma_{\uu}$, $\Omega\cup \Gamma_{\phi}$ are Gr\"oger-regular (see Definition \ref{def of geometric sets}).

\item \label{B2}$\C,\e,\vareps^0,\eps$ are essentially bounded on $\R^2$.
\item \label{B3}
There exists some $\alpha>0$ such that for all $\PP\in \R^2,\ \vareps\in \mathrm{Lin}_{\mathrm{sym}}(\R^2),\ \D\in \R^{2}$
\begin{subequations}
\begin{align*}
\C(\PP)\vareps:\vareps&\geq \alpha|\vareps|^2,\\
\eps(\PP)\D\cdot\D&\geq \alpha|\D|^2.\\
\end{align*}
\end{subequations}
\end{enumerate}
Then the Assumption \ref{A3} is satisfied.
\end{theorem}

\begin{remark}
By the setting of Theorem \ref{2d mixed}, a solution $(\uu,\phi,\PP)$ deduced from Theorem \ref{C^1 domain} satisfies
$$ (\uu(t),\phi(t),\PP(t))\in \big(W^{1,p}_{\Gamma_{\uu}}(\Omaga)\big)^d\times W^{1,p}_{\Gamma_{\phi}}(\Omaga)\times
\big(W^{1,p}_{\Gamma_{\PP}}(\Omaga)\big)^d$$
for a.e. $t\in(0,\hat{T})$. In particular, for a.e. $t\in(0,\hat{T})$ the solution $(\uu(t),\phi(t),\PP(t))$ has (componentwise) zero trace on the boundary portion $\Gamma_{\uu}\times\Gamma_{\phi}\times \Gamma_{\PP}$.
\end{remark}

Unfortunately, there is no such elliptic regularity result with $p>3$ by the 3D case when $(\uu,\phi)$ possesses mixed boundary conditions. In fact, to the best knowledge of the author, the sharpest value of $p$ for mixed case is exactly the one given by \cite{HDJonssonKneesRehberg2016}, which is merely expected to be close to two. Nevertheless, if one relaxes the conditions, in the sense that $(\uu,\phi)$ satisfies the Dirichlet boundary conditions, then an $L^p$-elliptic regularity result with $p>3$ is available (in fact, this is even true for all $p\in(2,\infty)$). Such $L^p$-regularity results can either be directly obtained using potential theory, or can also be obtained as a consequence of BMO-estimates for strongly elliptic systems. It is also worth noting that the Dirichlet $L^p$-elliptic regularity result is \textbf{regardless of the dimension}, meaning that it is applicable for both 2D and 3D cases. By our model, the $L^p$-elliptic regularity result is rephrased as follows:
\begin{theorem}\label{cube or C^1 domain}
Suppose that the following assumptions are satisfied:
\begin{enumerate}[label=(C\arabic*)]
\item \label{C1} $\Omaga\subset\R^d$ with $d\in\{2,3\}$ is an open rectangle (2D) or an open cuboid (3d) or a bounded domain with $C^1$-boundary, and $\Gamma_{\uu}=\Gamma_{\phi}=\partial\Omaga$.

\item \label{C2}$\C,\e,\vareps^0,\eps$ are continuous on $\R^d$ and there exists some continuous function $\alpha:\R^d\to (0,\infty)$ such that for all $\PP\in \R^d,\ \vareps\in \mathrm{Lin}_{\mathrm{sym}}(\R^d)$ and $\D\in \R^{d}$
\begin{subequations}
\begin{align*}
\C(\PP)\vareps:\vareps&\geq \alpha(\PP)|\vareps|^2,\\
\eps(\PP)\D\cdot\D&\geq \alpha(\PP)|\D|^2.\\
\end{align*}
\end{subequations}
\end{enumerate}
Then the Assumption \ref{A3} is satisfied.
\end{theorem}
The proof of Theorem \ref{cube or C^1 domain} will be a simple application of the regularity results from \cite{Acq1992,DolzmannMueller1995}. We refer to Section \ref{proof of thm cube or c1 domain} for details.
\begin{remark}
Compared to Theorem \ref{2d mixed}, we see that in Theorem \ref{cube or C^1 domain} no uniform boundedness and ellipticity conditions are imposed on the material parameters. This is due to the fact that the proof of the elliptic regularity result given in \cite{HDJonssonKneesRehberg2016} relies on a fixed point iteration, thus uniform upper and lower bounds of the coefficients will be essential for the proof therein. On the other hand, the proof of Theorem \ref{cube or C^1 domain} makes use of local perturbation arguments, we hence only ask for local boundedness and ellipticity.
\end{remark}
The particular interest is also devoted to 3D polyhedral domains, being one of the most widely used geometric objects in practical industrial applications. The elliptic regularity theory with symmetric coefficient tensor on polyhedral domains has a long standing history in mathematical research and has been nowadays rather completely developed. We refer the readers to \cite{MazyaRossmann2001,MazyaVladimirRossman2010} and the references therein for a self-contained study of the elliptic regularity theory on polyhedral domains. The main drawback of the methods given in \cite{MazyaRossmann2001,MazyaVladimirRossman2010} is that the elliptic regularity on a polyhedral domain is closely related to the spectral properties of the corresponding operator pencil, which is assumed to be self-adjoint, or in other words that the coefficient tensor is required to be symmetric. We recall that the coefficient tensor $\mathbb{B}(\PP)$ given by \eqref{not symmetric tentor} is in general asymmetric, which prevents us to apply the results from \cite{MazyaRossmann2001,MazyaVladimirRossman2010} to our model.

Nevertheless, the intrinsic physical principles of the ferroelectric model help us to solve this problem to some extent: a crucial observation from the actual ferroelectric model studied in \cite{GAMM2015} is that the material tensor $\mathbb{B}(\PP)$ satisfies the condition
\begin{align}\label{lame and laplace}
\mathbb{B}(\0)=\begin{pmatrix}\lambda {\mathbf{Id}_3}\otimes\mathbf{Id}_3+2\mu \mathbf{E}_3& \0 \\ \0& \gamma\mathbf{Id}_3\end{pmatrix},
\end{align}
where $\mathbf{Id}_3$ is the identity matrix in $\R^{3\times 3}$, $\mathbf{E}_3$ is the identity tensor in the space $\mathrm{Lin}_{\mathrm{sym}}(\R^{3\times 3})$ and $\lambda,\mu,\gamma$ are some given positive constants. A differential system with the upper left coefficient tensor in \eqref{lame and laplace} is the so called Lam\'{e} operator, while the bottom right coefficient tensor in \eqref{lame and laplace} corresponds to the Laplace operator. Physically, this phenomenon is to be understood that the piezo system decouples when the spontaneous polarization vanishes. We emphasize here that both Lam\'{e} and Laplace operators have \textbf{symmetric} coefficient tensors. The condition \eqref{lame and laplace} thus enables us to utilize the theory from \cite{MazyaVladimirRossman2010}, combining with a perturbation argument, to derive a reasonable $L^p$-elliptic regularity theory for our model. In order to formulate the $L^p$-elliptic regularity result for polyhedral domains we will still need the following assumption, which, roughly speaking, ensures that the results from \cite{MazyaVladimirRossman2010} are applicable for any ``singular point'' on the boundary of the polyhedral domain.
\begin{assumption}\label{sing}
Let $\Omaga\subset\R^3$ be a bounded polyhedral domain. Denote by $\mathrm{Sing}^1_{\Omaga}$ the set of singularities of a polyhedral domain $\Omaga$, i.e.
\begin{align*}
\mathrm{Sing}^1_{\Omaga}:=\{{\x}\in\partial\Omaga:{\x}\text{ is a vertex or lies on an edge}\}.
\end{align*}
We also define
\begin{align*}
\mathrm{Sing}^2_{\Omaga}:=\{{\x}\in\partial\Omaga:{\x}\text{ is a point on some open edge and $\theta_{{\x}}=\pi/2$}\},
\end{align*}
where $\theta_{{\x}}$ is the opening angle of the dihedron $\mathcal{D}_{{\x}}$ defined in Definition \ref{polyhedral domain}, namely, $\mathrm{Sing}^2_{\Omaga}$ is the set of edge points having $\frac{\pi}{2}$-opening angle;
\begin{multline*}
\mathrm{Sing}^3_{\Omaga}:=\{{\x}\in\partial\Omaga:{\x}\text{ is a vertex and }\\
\text{$\mathcal{K}_{{\x}}\cap B_1$ is $C^\infty$-diffeomorph to $\mathcal{K}_{c}\cap B_1$}\},
\end{multline*}
where $\mathcal{K}_{{\x}}$ is the cone defined in Definition \ref{polyhedral domain}, $B_1$ is the unit ball and $\mathcal{K}_{c}$ is the cone formed by the positive ${\x},{y},{z}$-axis, namely, $\mathrm{Sing}^3_{\Omaga}$ is the set of cube-corner-like vertices. Finally, we define
\begin{align*}
\mathrm{Sing}_{\Omaga}:=\mathrm{Sing}^1_{\Omaga}\setminus(\mathrm{Sing}^2_{\Omaga}\cup\mathrm{Sing}^3_{\Omaga}).
\end{align*}
Then we assume:
\begin{itemize}
\item For each vertex ${\x}$ the corresponding cone $\mathcal{K}_{{\x}}$ is a Lipschitz graph with edges, i.e. $\mathcal{K}_{{\x}}$ has the representation ${\x}_3>\psi({\x}_1,{\x}_2)$ in Cartesian coordinates with some Lipschitz function $\psi$.
\item $\mathrm{Sing}_{\Omaga}\subset \Gamma_{\PP}$.
\end{itemize}
\end{assumption}
\begin{remark}
Here follows an explanation why we need to impose additional regularity assumptions for vertices: One easily sees that every point on an open face or an open edge has an (in general small) neighborhood whose intersection with $\Omaga$ has Lipschitz boundary, and the validity of many important features such as Poincar\'e's inequality and Sobolev embedding etc. is guaranteed. However, a cone need not be a domain with Lipschitz boundary. For a counterexample, we refer to \cite[Fig. 9]{MazyaVladimirRossman2010}.
\end{remark}
We now state the last main result of this paper:
\begin{theorem}\label{poly thm}
Let the Assumption \ref{C2}, Assumption \ref{sing} and \eqref{lame and laplace} be satisfied. Let also $\Gamma_{\uu}=\Gamma_{\phi}=\partial\Omaga$. Then the Assumption \ref{A3} is satisfied.
\end{theorem}
\section{Notation and definitions}\label{notations}
\subsection{Basic notation}
The symbol $\mathrm{Lin}(A,B)$ ($\mathrm{Lin}_{\mathrm{sym}}(A)$) denotes the set of linear (linear symmetric) tensors from $A$ to $B$ (from $A$ to itself) for finite dimensional vector spaces $A$ and $B$. For matrices $M,N\in\R^{m\times n}$, the inner product $M:N$ of $M$ and $N$ is defined by $M:N:=\sum_{i=1}^m\sum_{j=1}^n M_{ij}N_{ij}$.

For a function ${f}:{\Lambda}\subset \R^m\rightarrow\R^n$, $D_{x}{f}(\cdot)\in\R^{n\times m}$ denotes the usual derivative of ${f}$ in the Euclidean space; for Banach spaces ${X},{Y}$ and a function $f:{X}\rightarrow {Y}$, the symbol $D_{{x}}f(z)[\bar{{x}}]$ denotes the G\^ateaux-differential of $f$ at point $z\in {X}$ in direction $\bar{{x}}\in {X}$. If ${\Lambda}$ is an open subset of $\R$, then the derivative of ${f}$ will also be denoted by ${f}'$, which stands for the time derivative of ${f}$. We also denote by $\nabla f$ the gradient of $f$, i.e. $\nabla f=(Df)^T$, where $(Df)^T$ is the transpose of the derivative of $f$.

For a differentiable matrix-valued function ${v}:\Omaga\subset\R^d\to \R^{d\times d}$, its divergence $\mathrm{div}\, {v}$ is defined as a $d$-dimensional vector with $(\mathrm{div}\,{v})_i:=\sum_{j=1}^d\partial_j{v}_{ij}$, i.e. the divergence is taken row-wise.

For Banach spaces ${X}$ and ${Y}$, we denote by $L({X},{Y})$ the space of all functions $f:{X}\to {Y}$ which are linear and continuous. $L({X},{Y})$ is a Banach space equipped with the norm $\|f\|_{L({X},{Y})}:=\sup_{\|x\|_{{X}}\leq 1}\|f(x)\|_{{Y}}$. We denote by $LH({X},{Y})$ the subset of $L({X},{Y})$ whose elements are additionally bijective.

We will also use the capital letter $C$ to denote a generic large positive number.
\subsection{Geometric settings}\label{geometric}
\subsubsection{$l$-set and Gr\"oger-regular set }\label{groeger sets}
In the following, we introduce the geometric concepts \textit{$l$-set} and \textit{Gr\"oger-regular set}, which plays a fundamental role in the proof of the main results:

\begin{definition}\label{def of geometric sets}
\begin{itemize}
\item[(i)]Let $d\in \N$ and $l\in(0,d]$. We say that a closed set $F\subset\R^d$ is an $l$-set, if for all $x\in F$ and $r\in(0,1]$ we have
    \begin{align*}
    \mathcal{H}_{l}(F\cap B(x,r))\sim r^{l},
    \end{align*}
where $\mathcal{H}_{l}$ denotes the $l$-dimensional Hausdorff measure.

\item[(ii)]
A set ${W}\subset\R^d$ is called Gr\"oger-regular, if $W$ is a bounded set and for every ${\x}\in\partial{W}$, there exist ${U}_1,{U}_2\subset\R^d$ and a bi-Lipschitz transformation $\Phi:{U}_1\rightarrow {U}_2$, such that ${\x}\in {U}_1$ and $\Phi({U}_1\cap W)$ is one of the following sets:
\begin{itemize}
\item[$\circ$] ${M}_1:=\{{\x}\in\R^d:|{\x}_i|<1\,\forall\, i=1,\cdots,d\text{ and }{\x}_{d}<0\}$,
\item[$\circ$] ${M}_2:=\{{\x}\in\R^d:|{\x}_i|<1\,\forall\, i=1,\cdots,d\text{ and }{\x}_{d}\leq 0\}$,
\item[$\circ$] ${M}_3:=\{{\x}\in {M}_2:{\x}_d<0\text{ or }{\x}_1> 0\}$.
\end{itemize}

\end{itemize}
\end{definition}

\begin{remark}
Here follow several comments on Definition \ref{def of geometric sets}:
\begin{itemize}
\item[(i)] The concept of $l$-set is also referred to as the \textit{$l$-Ahlfors-David-regular} set in literature, which was originally proved to be equivalent to the $L^2$-boundedness of the Cauchy integral operator appearing in the study of complex geometric analysis. We refer to the papers \cite{AhlforsPaper,DavidPaper,CauchyIntegralOperator} and the references therein for further details. After its first appearance, the Ahlfors-David-regular condition arises naturally in different research areas of mathematics and plays a fundamental role for many important results. We refer to \cite{jonssonWallin1984} for some enlightening examples of Ahlfors-David-regular-sets. Particularly, the closure of a bounded weakly Lipschitz domain (in the sense of \cite[Def. 1.2.1.2]{grisvard2011}) in $\R^d$ is a $d$-set and the boundary of a bounded Lipschitz domain in $\R^d$ is a $(d-1)$-set, see for instance \cite[Chap. II, Exp. 1]{jonssonWallin1984} for a proof.

\item[(ii)] The Gr\"oger-regular set was originally introduced by Gr\"oger \cite{GroegerRehberg1989} for the purpose of characterizing the geometric profile of domains of mixed type. It turns out that the Gr\"oger-regular condition is not a too restrictive condition and  therefore widely applied in different applications. We refer to \cite[Sec. 7]{Haller_Dintelmann_2009} for a survey on several representative examples of Gr\"oger-regular sets.

\end{itemize}
\end{remark}

A schematic description for Gr\"oger-regular set is given by Fig. \ref{groeger} below.
\begin{figure}[!htbp]
  \centering
  \includegraphics[width=50mm]{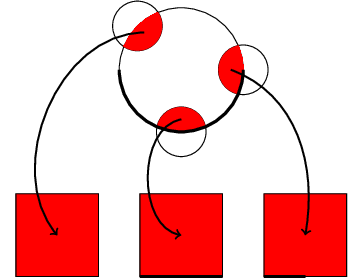}
  \caption{A schematic description of a Gr\"oger-regular set. From left to right: the sets $M_1$, $M_2$ and $M_3$.}\label{groeger}
\end{figure}

\subsubsection{Polyhedral domains}
We follow \cite{MazyaVladimirRossman2010} to define a polyhedral domain. The definitions of a dihedron and a cone will also be given, which form part of the definition of a polyhedral domain.
\begin{definition}[Dihedron]
Let $(r,\phi)\in(0,\infty)\times (0,2\pi]$ be the polar coordinate of a non-zero two dimensional point ${\x}'=({\x}_1,{\x}_2)$, i.e. $r=|{\x}'|$ and ${\x}'=(r\cos \phi,r\sin\phi)^T$. The two dimensional wedge $\mathcal{W}$ with opening angle $\theta\in(0,2\pi]$ is defined by
\begin{align*}
\mathcal{W}:=\{{\x}'=({\x}_1,{\x}_2):r\in(0,\infty),\phi\in(0,\theta)\}.
\end{align*}
The dihedron $\mathcal{D}\subset\R^3$ with opening angle $\theta$ is then defined by $\mathcal{D}:=\mathcal{W}\times \R$.
\end{definition}

\begin{definition}[Cone]
A cone $\mathcal{K}\subset\R^3$ is defined by
\begin{align*}
\mathcal{K}:=\{{\x}\in\R^3-\{\0\}:{\x}/|{\x}|\in{\Gamma}\},
\end{align*}
where ${\Gamma}$ is a subdomain on the unit sphere $\mathbb{S}^{2}$ such that
\begin{align*}
\partial{\Gamma}=\bar{{\Gamma}}_1\cup...\cup\bar{{\Gamma}}_k
\end{align*}
for some $k\in\N_{\geq 3}$, where ${\Gamma}_i$ are pairwise disjoint and non-curvewise collinear open arcs (namely, the intersection point of two different arcs has locally different tangents along the arcs) on the unit sphere.
\end{definition}
\begin{definition}[Polyhedral domain]\label{polyhedral domain}
A bounded domain $\Omaga\subset\R^3$ is said to be a polyhedral domain, if
\begin{itemize}
\item The boundary $\partial\Omaga$ is a disjoint union of smooth open two dimensional manifolds $\mathcal{F}_j,j=1,...,l$ (the faces of $\Omaga$), smooth open curves $K_{j},j=1,...,m$ (the edges of $\Omaga$) and vertices ${\x}^j\in\R^3,j=1,...,k$.
\item For every ${\x}\in K_{j}$ there exist a neighborhood ${U}_{{\x}}\subset\R^3$ of ${\x}$ and a $C^\infty$-diffeomorphism ${\iota}_{{\x}}$ such that ${\iota}_{{\x}}$ maps ${U}_{{\x}}\cap \Omaga$ onto $\mathcal{D}_{{\x}}\cap B_1$, where $\mathcal{D}_{{\x}}$ is a dihedron and $B_1$ is the unit ball in $\R^3$.
\item For every vertex ${\x}^j$ there exist a neighborhood ${U}_{j}\subset\R^3$ of ${\x}^j$ and a $C^\infty$-diffeomorphism ${\iota}_{j}$ such that ${\iota}_j$ maps ${U}_{j}\cap\Omaga$ onto $\mathcal{K}_{j}\cap B_1$, where $\mathcal{K}_j$ is a cone.
\end{itemize}
\end{definition}

\subsection{Function spaces}
The Sobolev space $W^{k,p}(\Omaga)$ with non-negative integer $k$ and real number $p\in[1,\infty]$ is defined by
\begin{equation*}
W^{k,p}(\Omaga):=\{f\in L^p(\Omaga):D^{{\beta}} f\in L^p(\Omaga)\ \mathrm{for\ }|{\beta}|\leq k\}
\end{equation*}
with the norm \begin{equation*}
\|f\|_{W^{k,p}(\Omaga)}:=\big(\sum_{|{\beta}|\leq k}\|D^{{\beta}} f\|^p_{L^p(\Omaga)}\big)^{\frac{1}{p}},
\end{equation*}
where $D^{{\beta}}$ denotes the multi-differential symbol for a non-negative integer multi-index ${{\beta}}$. For a closed set ${\Gamma}\subset\partial\Omaga$, the space $W^{1,p}_{{\Gamma}}(\Omaga)$ with $p\in(1,\infty)$ is defined as the closure of the space
\begin{equation*}
\{\phi|_{\Omaga}:\phi\in C_0^\infty(\R^d),\text{supp}(\phi)\cap {\Gamma}=\varnothing\}
\end{equation*}
in the space $W^{1,p}(\Omaga)$ w.r.t. the $W^{1,p}$-norm. The dual space of $W^{1,p}_{{\Gamma}}(\Omaga)$ is denoted by $W^{-1,p'}_{{\Gamma}}(\Omaga)$, where $p'$ is the H\"older conjugate of $p$ with $\frac{1}{p}+\frac{1}{p'}=1$. We also write $H^1_{{\Gamma}}(\Omaga)=W^{1,2}_{{\Gamma}}(\Omaga)$ and its dual is denoted by $H^{-1}_{{\Gamma}}(\Omaga)$. For ${\Gamma}=\partial\Omaga$, $W^{1,p}_{{\Gamma}}(\Omaga)$ and $H^1_{{\Gamma}}(\Omaga)$ are denoted by $W^{1,p}_{0}(\Omaga)$ and $H^1_{0}(\Omaga)$ and their dual spaces are denoted by $W^{-1,p'}(\Omaga)$ and $H^{-1}(\Omaga)$ respectively. For ${f}\in \big(L^p(\Omaga)\big)^d$ the distributional divergence operator $\mathrm{Div}:\big(L^p(\Omaga)\big)^d\rightarrow W_{\Gamma}^{-1,p}(\Omaga)$ is defined by
\begin{align*}
\mathrm{Div}\,{f}[g]:=-\int_{\Omaga}{f}\cdot \nabla g\,d{\x}
\end{align*}
for all $g\in W_{\Gamma}^{1,p'}(\Omaga)$. Throughout the paper, we will also use the following shorthand notation of function spaces:
\begin{align}
\W&:=\big(W^{1,p}_{\Gamma_{\uu}}(\Omaga)\big)^d\times W^{1,p}_{\Gamma_{\phi}}(\Omaga),\label{special not1}\\
\W^{*}&:=\big(W^{-1,p}_{\Gamma_{\uu}}(\Omaga)\big)^d\times W^{-1,p}_{\Gamma_{\phi}}(\Omaga),\label{special not2}\\
\V&:=(\big(W^{1,p}_{\Gamma_{\PP}}(\Omaga)\big)^d,\big(W^{-1,p}_{\Gamma_{\PP}}(\Omaga)\big)^d)_{\frac{1}{r},r},\label{special not3}\\
\U&:=\{\,\PP\in \big(C(\yoverline{\Omaga})\big)^d:\PP|_{\Gamma_{\PP}}=\0\,\}.\label{special not4}
\end{align}

\subsection{Real and complex interpolation}
Let $\theta,s\in(0,1),q\in [1,\infty]$ and let ${X},{Y}$ be compatible Banach spaces in the sense of \cite[Chap. 2.3]{Bergh1976}. Then $({X},{Y})_{\theta,q}$ and $[{X},{Y}]_s$ denote the real and complex interpolation space induced by $({X},{Y})$ of components $(\theta,q)$ and of component $s$ respectively. Due to their cumbersome construction we refer to \cite{Bergh1976,Triebel1978} for the definitions and properties of real and complex interpolation spaces.

\section{Proof of the local well-posedness result}\label{proof of thm 1}
\subsection{An abstract existence and uniqueness result for quaislinear parabolic equations}
Our local well-posedness result relies on the existence and uniqueness result given by \cite{LiClement1994}. To introduce the result, we need the following definition:
\begin{definition}[Maximal parabolic regularity]\label{maximal parabolic property}
Let $J=(T_0,T_1)$ be a time interval. Let $X$ be a Banach space and $A:\mathrm{dom}(A)\to X$ be a closed linear operator with dense domain $\mathrm{dom}(A)\subset X$. Suppose ${\tau} \in (1,\infty)$. Then we say that $A$ has maximal parabolic $L^{\tau} (J;X)$-regularity if and only if for every $f\in L^{\tau} (J;X)$ there is a unique function $\mathbf{w}\in W ^{1,{\tau}} (J;X) \cap L^{\tau} \big(J;\mathrm{dom}(A)\big)$ which satisfies $\mathbf{w}(T_0 ) = \0$ and
\begin{alignat*}{2}
&\mathbf{w} '(t) + A\mathbf{w}(t) = f(t)\quad&&\mathrm{in\ }X
\end{alignat*}
for a.a. $t \in J$.
\end{definition}
\begin{remark}
It is shown in \cite{Dore1993} that if $A$ has maximal parabolic $L^{{\tau}_0} (J;X)$-regularity for some ${\tau}_0\in(1,\infty)$, then it has maximal parabolic $L^{\tau} (J;X)$-regu\-larity for all ${\tau}\in(1,\infty)$. It is therefore unambiguous to merely say that $A$ has maximal parabolic regularity.
\end{remark}
The maximal parabolic regularity property motivates the following existence and unique result introduced in \cite{LiClement1994} (see also \cite{pruss2002}) for a quasilinear parabolic system, which plays the main role in the analysis of our local well-posedness result:
\begin{theorem}\label{pruss2002}
Let $Y,X$ be Banach spaces, $Y\hookrightarrow X$ densely and let $\tau \in (1,\infty)$. Suppose that $A:J\times (Y,X)_{\frac{1}{\tau},\tau}\to L(Y ,X)$ is continuous and $A(0,\mathbf{w}_0 )$ satisfies maximal parabolic regularity on $X$ with $\mathrm{dom}\big(A(0,\mathbf{w}_0 )\big) = Y$ for some $\mathbf{w}_0 \in (Y,X)_{\frac{1}{\tau},\tau}$. Let $S: J \times (Y,X)_{\frac{1}{\tau},\tau}\rightarrow X$
be a Carath\'{e}odory map, i.e., $S(\cdot,x)$ is Bochner-measurable for each $x\in (Y,X)_{\frac{1}{\tau},\tau}$ and $S(t,\cdot)$ is continuous for a.a. $t\in J$. Moreover, let $S(\cdot,\0)$ be from $L^\tau (J;X)$ and the following assumptions be satisfied:
\begin{itemize}
\item[(A)] For every $M > 0$, there exists a positive constant $C_M$ such that for all $t\in J$ and all $\mathbf{w}, \bar{\mathbf{w}}\in (Y,X)_{\frac{1}{\tau},\tau}$ with $\max(\|\mathbf{w}\|_{(Y,X)_{\frac{1}{\tau},\tau}},\|\bar{\mathbf{w}}\|_{(Y,X)_{\frac{1}{\tau},\tau}}) \leq M$, we have
\begin{align*}
\|A(t,\mathbf{w})-A(t,\bar{\mathbf{w}})\|_{L(Y,X)}\leq C_M\|\mathbf{w}-\bar{\mathbf{w}}\|_{(Y,X)_{\frac{1}{\tau},\tau}};
\end{align*}
\item[(S)] For every $M > 0$ there exists a function $h_M \in L^\tau (J)$ such that for all $\mathbf{w}, \bar{\mathbf{w}}\in (Y,X)_{\frac{1}{\tau},\tau}$ with $\max(\|\mathbf{w}\|_{(Y,X)_{\frac{1}{\tau},\tau}},\|\bar{\mathbf{w}}\|_{(Y,X)_{\frac{1}{\tau},\tau}}) \leq M$, it is true that
\begin{align*}
\|S(t,\mathbf{w})-S(t,\bar{\mathbf{w}})\|_{X}\leq h_M(t)\|\mathbf{w}-\bar{\mathbf{w}}\|_{(Y,X)_{\frac{1}{\tau},\tau}}
\end{align*}
for a.a. $t \in J$.
\end{itemize}
Then there exists some maximal $T_{\mathrm{max}}\in J$ such that the problem
\begin{alignat*}{2}
\mathbf{w}' (t) + A\big(t,\mathbf{w}(t)\big)\mathbf{w}(t)& = S\big(t,\mathbf{w}(t)\big) \quad&&\mathrm{in\ } J\times X
\end{alignat*}
with $\mathbf{w}(T_0 ) = \mathbf{w}_0$ admits a unique solution
\begin{align}\label{maximality of hat T}
\mathbf{w}\in W^{1,\tau} (0 ,S ;X) \cap L^\tau (0 ,S;Y )
\end{align}
on $(0,S)$ for every $S\in (0,T_{\mathrm{max}})$. The maximality of $T_{\mathrm{max}}$ is understood as follows:
\begin{itemize}
\item Either \eqref{maximality of hat T} does not hold for $S=T_{\mathrm{max}}$, or
\item $T_{\mathrm{max}}=T_1$ and \eqref{maximality of hat T} also holds for $S=T_1$.
\end{itemize}
\end{theorem}
\subsection{Embedding results for interpolation spaces}
The interpolation space $\V$ defined by \eqref{special not3} displays a somehow complicated structure at the first glance, due to the cumbersome definition of interpolation spaces. We show that $\V$ is continuously embedded into some H\"older space, which will simplify several arguments in the proof of the main results.
\begin{lemma}\label{linfty}
Let the Assumption \ref{A1} be satisfied. Then for all $p\in(d,\infty)$ and $r\in(\frac{2p}{p-d},\infty)$ we have the continuous embedding
\begin{align}
\V\hookrightarrow \big(W^{1-2\mathfrak{p},p}_{\Gamma_{\PP}}(\Omaga)\big)^d
\end{align}
for some $\mathfrak{p}\in(\frac{1}{r},\frac{p-d}{2p})$. Consequently, we have the embedding
\begin{align}\label{tau embedding}
\big(W^{1-2\mathfrak{p},p}_{\Gamma_{\PP}}(\Omaga)\big)^d\hookrightarrow \big(C^\delta(\yoverline{\Omaga})\big)^d
\end{align}
with $\delta=1-2\mathfrak{p}-\frac{d}{p} \in(0,1)$.
\end{lemma}
\begin{proof}
The statement is originally proved in \cite[Lem. A.1]{MeinschmidtMeyerRehberg16} for the case $d=3$ and $\Gamma_{\PP}=\varnothing$. We adopt the idea therein to prove the statement for our model. From the condition $\mathfrak{p}\in(\frac{1}{r},\frac{p-d}{2p})$ we obtain that
$$\delta=1-2\mathfrak{p}-\frac{d}{p} \in(0,1).$$
Thus \eqref{tau embedding} follows immediately from Rellich's embedding theorem. It is left to show that there exists some $\mathfrak{p}\in (\frac{1}{r},\frac{p-d}{2p})$ such that
$$\V\hookrightarrow\big(W^{1-2\mathfrak{p},p}_{\Gamma_{\PP}}(\Omaga)\big)^d.  $$
Let $\nu\in(0,\frac{p-d}{2p})$. We obtain that
\begin{align}\label{azero}
&(\big(W_{\Gamma_{\PP}}^{1,p}(\Omaga)\big)^d,\big(W_{\Gamma_{\PP}}^{-1,p}(\Omaga)\big)^d)_{\nu,1}\nonumber\\
=&(\big(W_{\Gamma_{\PP}}^{1,p}(\Omaga)\big)^d,\Big(\big(W_{\Gamma_{\PP}}^{1,p}(\Omaga)\big)^d,\big(W_{\Gamma_{\PP}}^{-1,p}(\Omaga)\big)^d\Big)_{\frac{1}{2},1})_{2\nu,1}\nonumber\\
=&(\big(W_{\Gamma_{\PP}}^{1,p}(\Omaga)\big)^d,\Big(\big(W_{\Gamma_{\PP}}^{-1,p}(\Omaga)\big)^d,\big(W_{\Gamma_{\PP}}^{1,p}(\Omaga)\big)^d\Big)_{\frac{1}{2},1})_{2\nu,1}\nonumber\\
=&(\big(W_{\Gamma_{\PP}}^{1,p}(\Omaga)\big)^d,\Big(\big(W_{\Gamma_{\PP}}^{-1,p}(\Omaga)\big)^d,\mathrm{dom}(-\Delta+1)\Big)_{\frac{1}{2},1})_{2\nu,1},
\end{align}
where the first equality comes from the reiteration theorem given in  \cite[Chap. 1.10.2]{Triebel1978}, the second from the property
\begin{align}\label{thetaq}
(X,Y)_{\theta,q}=(Y,X)_{1-\theta,q}
\end{align}
of real interpolation spaces, and the last equality from the fact that
\begin{equation*}
\big(W_{\Gamma_{\PP}}^{1,p}(\Omaga)\big)^d=\mathrm{dom}(-\Delta+1),
\end{equation*}
where $\mathrm{dom}(-\Delta+1)$ is defined as the domain of the operator
$$-\Delta+1:\big(W_{\Gamma_{\PP}}^{1,p}(\Omaga)\big)^d\to\big(W_{\Gamma_{\PP}}^{-1,p}(\Omaga)\big)^d.$$
According to \cite[Thm. 11.5]{SquareRoot}, the extension of the operator $-\Delta+1$ from $\big(L^p(\Omaga)\big)^d$ to $\big(W_{\Gamma_{\PP}}^{-1,p}(\Omaga)\big)^d$ is a positive operator in the sense of \cite[Chap. 1.14]{Triebel1978}. This property enables us to exploit \cite[Chap. 1.15.2]{Triebel1978} to define the square root $(-\Delta+1)^{\frac{1}{2}}$ of $-\Delta+1$ on $\big(W_{\Gamma_{\PP}}^{-1,p}(\Omaga)\big)^d$ with domain $\mathrm{dom}\big((-\Delta+1)^{\frac{1}{2}}\big)$ and to obtain that
\begin{align}\label{aoneone}
\Big(\big(W_{\Gamma_{\PP}}^{-1,p}(\Omaga)\big)^d,\mathrm{dom}(-\Delta+1)\Big)_{\frac{1}{2},1}\hookrightarrow\mathrm{dom}\big((-\Delta+1)^{\frac{1}{2}}\big).
\end{align}
From \cite[Thm. 5.1 (ii)]{SquareRoot} one infers the embedding
\begin{align}\label{athree}
\mathrm{dom}\big((-\Delta+1)^{\frac{1}{2}}\big)\hookrightarrow \big(L^{p}(\Omaga)\big)^d.
\end{align}
Finally, from \eqref{aoneone} and \eqref{athree} we conclude that
\begin{align}\label{afour}
\Big(\big(W_{\Gamma_{\PP}}^{-1,p}(\Omaga)\big)^d,\mathrm{dom}(-\Delta+1)\Big)_{\frac{1}{2},1}\hookrightarrow\big(L^{p}(\Omaga)\big)^d.
\end{align}
It follows that
\begin{align}
&(\big(W_{\Gamma_{\PP}}^{1,p}(\Omaga)\big)^d,\big(W_{\Gamma_{\PP}}^{-1,p}(\Omaga)\big)^d)_{\nu,1}\nonumber\\
\hookrightarrow&(\big(W_{\Gamma_{\PP}}^{1,p}(\Omaga)\big)^d,\big(L^{p}(\Omaga)\big)^d)_{2\nu,1}\nonumber\\
=&(\big(L^{p}(\Omaga)\big)^d,\big(W_{\Gamma_{\PP}}^{1,p}(\Omaga)\big)^d)_{1-2\nu,1}\nonumber\\
\hookrightarrow&[\big(L^{p}(\Omaga)\big)^d,\big(W_{\Gamma_{\PP}}^{1,p}(\Omaga)\big)^d]_{1-2\nu}=\big(W_{\Gamma_{\PP}}^{1-2\nu,p}(\Omaga)\big)^d,\label{koukou}
\end{align}
where the first embedding comes from \eqref{azero} and \eqref{afour}, the first equality is obtained by using \eqref{thetaq}, the second embedding comes from the property
\begin{align*}
(X,Y)_{\theta,1}\hookrightarrow [X,Y]_{\theta}
\end{align*}
of interpolation spaces, and the second equality is deduced from \cite[Thm. 3.1]{Griepentrog_2002} (notice that $\nu$ is in $(0,\frac{p-d}{2p})$, which particularly implies that $1-2\nu>\frac{d}{p}>\frac{1}{p}$ and consequently that $1-2\nu\neq\frac{1}{p}$, being the condition of \cite[Thm. 3.1]{Griepentrog_2002}). Finally, using the fundamental properties of real interpolation spaces we obtain that
\begin{align*}
&(\big(W_{\Gamma_{\PP}}^{1,p}(\Omaga)\big)^d,\big(W_{\Gamma_{\PP}}^{-1,p}(\Omaga)\big)^d)_{\frac{1}{r},r}\\
\hookrightarrow&(\big(W_{\Gamma_{\PP}}^{1,p}(\Omaga)\big)^d,\big(W_{\Gamma_{\PP}}^{-1,p}(\Omaga)\big)^d)_{\nu,1}\\
\hookrightarrow& \big(W_{\Gamma_{\PP}}^{1-2\nu,p}(\Omaga)\big)^d
\end{align*}
for $\nu\in(\frac{1}{r},1)$, where the first embedding follows from \cite[1.3.3. (4)]{Triebel1978}. The condition on $r$ implies that the interval $(\frac{1}{r},\frac{p-d}{2p})$ is not empty, thus one can choose some $\mathfrak{p}\in(\frac{1}{r},\frac{p-d}{2p})\subset(0,\frac{p-d}{2p})$ such that \eqref{koukou} is valid. This completes the proof.
\end{proof}

The following proposition reveals the temporal and spatial H\"older continuity of a local solution $\PP$.
\begin{proposition}\label{reg of polar}
Let the Assumption \ref{A1} be satisfied. Then for $p\in(d,\infty)$ and $r\in(\frac{2p}{p-d},\infty)$ there exists some $\delta\in(0,1)$ such that the embedding
\begin{align*}
W^{1,r}\Big(J;\big(W^{-1,p}_{\Gamma_{\PP}}(\Omaga)\big)^d\Big)\cap L^{r}\Big(J;\big(W^{1,p}_{\Gamma_{\PP}}(\Omaga)\big)^d\Big)\hookrightarrow C^{\delta}\Big(\bar{J};\big(C^{\delta}(\yoverline{\Omaga})\big)^d\Big)
\end{align*}
is valid for each interval $J\subset\R$.
\end{proposition}
\begin{proof}
From \cite[Thm. 3]{amann2001} we have the embedding
\begin{align*}
&W^{1,r}\Big(J;\big(W^{-1,p}_{\Gamma_{\PP}}(\Omaga)\big)^d\Big)\cap L^{r}\Big(J;\big(W^{1,p}_{\Gamma_{\PP}}(\Omaga)\big)^d\Big)\nonumber\\
\hookrightarrow&\, C^{s-\frac{1}{r}}\Big(\bar{J};\Big(\big(W^{-1,p}_{\Gamma_{\PP}}(\Omaga),W^{1,p}_{\Gamma_{\PP}}(\Omaga)\big)_{\theta,1}\Big)^d\Big)
\end{align*}
for $s\in(\frac{1}{r},1)$ and $\theta\in[0,1-s)$. Notice that $1-\theta\in(s,1]$ and $s\in(\frac{1}{r},1)$ and $1-\theta,s$ can be arbitrarily chosen in these intervals. The condition on $r$ implies that $\frac{1}{r}<\frac{p-d}{2p}$, thus one can choose $s$ sufficiently small and then $\theta$ sufficiently large such that $s\in(\frac{1}{r},\frac{p-d}{2p})$ and $1-\theta\in(s,\frac{p-d}{2p}]\subset(0,\frac{p-d}{2p})$. Then \eqref{koukou} is satisfied and we obtain that
\begin{align*}
\big(W^{-1,p}_{\Gamma_{\PP}}(\Omaga),W^{1,p}_{\Gamma_{\PP}}(\Omaga)\big)_{\theta,1}&=\big(W^{1,p}_{\Gamma_{\PP}}(\Omaga),W^{-1,p}_{\Gamma_{\PP}}(\Omaga)\big)_{1-\theta,1}\\
&\hookrightarrow C^{1-2(1-\theta)-\frac{d}{p}}(\yoverline{\Omaga}).
\end{align*}
By choosing $\delta:=\min\{s-\frac{1}{r},1-2(1-\theta)-\frac{d}{p}\}\in(0,1)$ we obtain the desired result.
\end{proof}
\subsection{Uniform boundedness of \protect{${\LL}^{-1}_{\mathbf{P}}$}}
Let the number $p$ and the space $\U$ be as defined by Assumption \ref{A3} and let $(\uu,\phi)$ be the solution of the differential equation
\begin{align*}
\LL_{\PP}(\uu,\phi)=\lol
\end{align*}
for $\PP\in\U$ and $\lol\in \W^*$. Then from Assumption \ref{A3} we have the estimate
\begin{align}\label{haochang4}
\|(\uu,\phi)\|_{\W}\leq C_{\PP}\|\lol\|_{\W^*},
\end{align}
where $C_{\PP}=\|\LL^{-1}_{\PP}\|_{L(\W^*,\W)}$ is the inverse operator norm. Under the Assumptions \ref{B1} to \ref{B3} the constant $C_{\PP}$ is uniform for all test functions $\PP\in \U $, which is a direct consequence of Theorem \ref{higherorder}. In general, however, the inverse operator norm $C_{\PP}$ also depends on the modulus of continuity of $\mathbb{B}(\PP)$, thus $C_{\PP}$ could tend to infinity by testing various $\PP$ even when they are uniformly bounded in $\U$, which makes the verification of the condition (S) of Theorem \ref{pruss2002} impossible. From Lemma \ref{linfty} we know yet that all test functions $\PP$ will be taken from a bounded subset of a H\"older space. Using the compact embedding of a H\"older space into $\big(C(\yoverline{\Omaga})\big)^d$ and the indirect continuity arguments given by \cite{MeinschmidtMeyerRehberg16} we are able to show the uniform boundedness of $C_{\PP}$.
\begin{lemma}\label{haochang6}
Let the Assumptions \ref{A1} to \ref{A3} be satisfied and let $p$ be the number given by Assumption \ref{A3}. Let $\mathcal{M}$ be a compact subset of the space $\U$. Then the norm $C_{\PP}$ of the inverse operator $\mathrm{L}^{-1}_{\PP}$ is uniformly bounded by some positive constant $C^*$ for all $\PP\in \mathcal{M}$.
\end{lemma}
\begin{proof}
Firstly we show that the mapping
\begin{align}\label{continuous mapping}
\U\ni\PP\mapsto J(\PP):=\left(-\mathrm{Div}\left(\mathbb{B}(\PP)\begin{pmatrix}\vareps(\cdot|_{\uu}) \\ \nabla\cdot|_{\phi}\end{pmatrix}\right)\right)^{-1}\in LH\big(\W^*,\W\big)
\end{align}
is well-defined and continuous, where $LH(X,Y)$ denotes the set of linear homeomorphisms between Banach spaces $X$ and $Y$. On the one hand, we obtain that
\begin{align}\label{pp1-pp2 6}
&\mathrm{L}_{\PP_1}(\uu,\phi)[(\buu,\bphi)]-\mathrm{L}_{\PP_2}(\uu,\phi)[(\buu,\bphi)]\nb\\[4pt]
=&\int_{\Omaga}(\mathbb{B}(\PP_1)-\mathbb{B}(\PP_2))\begin{pmatrix}\vareps(\uu)\\ \nabla\phi\end{pmatrix}:\begin{pmatrix}\vareps(\buu)\\ \nabla\bphi\end{pmatrix}\,dx\nb\\[4pt]
\leq&\, C^*_1\int_{\Omaga}|\PP_1-\PP_2|\cdot|(\vareps(\uu),\nabla\phi)|\cdot|(\vareps(\buu),\nabla\bphi)|\,dx\nb\\[4pt]
\leq&\,C C^*_1\|\PP_1-\PP_2\|_{L^\infty(\Omaga)}\|(\uu,\phi)\|_{\W}\|(\buu,\bphi)\|_{\WW}
\end{align}
for all $(\uu,\phi)\in \W$ and $(\buu,\bphi)\in \WW$, where $C^*_1$ denotes the Lipschitz constant of $\mathbb{B}$ on the closed ball $\bar{B}_M(0)\subset\R^d$, and $M$ is the upper bound of $\mathcal{M}$ in $\U$. From this we infer that the mapping
\begin{align*}
\U\ni\PP\mapsto-\mathrm{Div}\left(\mathbb{B}(\PP)\begin{pmatrix}\vareps(\cdot|_{\uu}) \\ \nabla\cdot|_{\phi}\end{pmatrix}\right)\in LH\big(\W,\W^*\big)\nonumber
\end{align*}
is continuous (that the image is in fact a linear homeomorphism follows from Assumption \ref{A3}). On the other hand, from fundamental analysis we know that for Banach spaces $X$ and $Y$, the mapping $LH(X,Y)\ni B\mapsto B^{-1}\in LH(Y,X)$ is continuous. Then the continuity of \eqref{continuous mapping} follows, as the mapping $J$ is a composition of continuous functions. Now since $J$ is continuous and $\mathcal{M}$ is compact, we see that $J(\mathcal{M})$ is a compact subset of $L\big(\W^*,\W\big)$. In particular we obtain that
$$C^*=\sup_{\PP\in\mathcal{M}}\|J(\PP)\|_{L(\W^*,\W)}<\infty. $$
This completes the proof.
\end{proof}
\subsection{Proof of Theorem \ref{C^1 domain}}
We are now ready to prove Theorem \ref{C^1 domain}:
\begin{proof}[Proof of Theorem \ref{C^1 domain}]
Let us first recall the following function spaces
\begin{align*}
\W&=\big(W^{1,p}_{\Gamma_{\uu}}(\Omaga)\big)^d\times W^{1,p}_{\Gamma_{\phi}}(\Omaga),\\
\W^{*}&=\big(W^{-1,p}_{\Gamma_{\uu}}(\Omaga)\big)^d\times W^{-1,p}_{\Gamma_{\phi}}(\Omaga),\\
\V&=(\big(W^{1,p}_{\Gamma_{\PP}}(\Omaga)\big)^d,\big(W^{-1,p}_{\Gamma_{\PP}}(\Omaga)\big)^d)_{\frac{1}{r},r},\\
\U&=\{\,\bar{\PP}\in \big(C(\yoverline{\Omaga})\big)^d:\bar{\PP}|_{\Gamma_{\PP}}=\0\,\}
\end{align*}
defined by \eqref{special not1} to \eqref{special not4}, which will be extensively used in the remaining proof. We formulate the following notation corresponding to the ones given by Theorem \ref{pruss2002}:
\begin{align*}
A(t,\PP)&\equiv-\Delta,\\
\tau&=r,\\
Y&=\big(W_{\Gamma_{\PP}}^{1,p}(\Omaga)\big)^d, \\
X&=\big(W_{\Gamma_{\PP}}^{-1,p}(\Omaga)\big)^d
\end{align*}
and
\begin{align}\label{stp}
S(t,\PP)=\sos \big(t,\vareps(\uu_{\PP}(t)),\nabla\phi_{\PP}(t),\PP \big)
\end{align}
where $\big(\uu_{\PP}(t),\phi_{\PP}(t)\big)$ is the unique $\W$-weak solution of the differential equation
\begin{align}\label{uphi}
{\LL_{\PP}}\big(\uu_{\PP}(t),\phi_{\PP}(t)\big)={\ell (}t,\PP)
\end{align}
for $\PP\in(Y,X)_{\frac{1}{\tau},\tau}=\V\subset\U$ (the inclusion is deduced from Lemma \ref{linfty}), with ${\ell (}t,\PP)$ and $\sos \big(t,\vareps,\nabla\phi,\PP \big)$ given by \eqref{ltp}, \eqref{stp intro} respectively. In particular, $\big(\uu_{\PP}(t),\phi_{\PP}(t)\big)$ is uniquely determined by a given pair $(t,\PP)$, thus $S(t,\PP)$ is well-defined. Having defined these notation, we utilize Theorem \ref{pruss2002} to show that the equation
$$\PP'(t)+A\big(t,\PP(t)\big)=S\big(t,\PP(t)\big) $$
with initial value $\PP_0$ has a unique local solution $\PP$. We firstly give the following statements corresponding to part of the conditions from Theorem \ref{pruss2002}, which are relatively easier to verify:
\begin{itemize}
\item That $Y\hookrightarrow X$ densely is trivial.
\item That $A(0,\PP_0)=-\Delta$ satisfies maximal parabolic regularity on $X$ with $\mathrm{dom}\big(A(0,\PP_0)\big)=Y$ follows from \cite[Thm. 5.4, Rem. 5.5]{Haller_Dintelmann_2009}, where by applying \cite[Thm. 1.2]{Egert2018} the volume-preserving condition in \cite{Haller_Dintelmann_2009} can be weakened to the $(d-1)$-Ahlfors-David-regular condition.

\item That $A(t,\PP)$ is continuous from $[0,T]\times (Y,X)_{\frac{1}{\tau},\tau}$ to $L(Y,X)$ is trivial, since $A$ is constantly equal to $-\Delta$.
\item Since $A$ is constantly equal to $-\Delta$, the validity of (A) in Theorem \ref{pruss2002} is evident.
\item Now we show that
\begin{align*}
S(t,\cdot):(Y,X)_{\frac{1}{\tau},\tau}\to X
\end{align*}
is continuous for a.a. $t\in[0,T]$. Let $\PP\in (Y,X)_{\frac{1}{\tau},\tau}$ be arbitrary. For $\kappa>0$ define
\begin{equation*}
\mathcal{O}_{\kappa}(\PP):=\{\BPP\in (Y,X)_{\frac{1}{\tau},\tau}:\|\PP-\BPP\|_{(Y,X)_{\frac{1}{\tau},\tau}}<\kappa\}.
\end{equation*}
Thus to show the continuity, it suffices to show that
\begin{align*}
\lim_{\kappa\to 0}\sup_{\BPP\in \mathcal{O}_{\kappa}(\PP)}\|S(t,\PP)-S(t,\BPP)\|_{X}=0.
\end{align*}
But this follows directly from the Assumption (S) given in Theorem \ref{pruss2002}, which will be shown in the subsequent part of the proof below and we do not repeat here.
\end{itemize}
It is left to show that
\begin{itemize}
\item given a fixed $\PP$ in $(Y,X)_{\frac{1}{\tau},\tau}$, the mapping $t\mapsto S(t,\PP)$ is Bochner-measur\-able,
\item the validity of (S) in Theorem \ref{pruss2002} and
\item $S(\cdot,\0)$ is from $L^r(0,T;X)$.
\end{itemize}
We will show these statements in the following steps.
\subsubsection*{Step 1: Bochner-measurability of $t\mapsto S(t,\PP)$}
Recall that
\begin{align*}
S(t,\PP)=&\,\sos \big(t,\vareps(\uu_{\PP}(t)),\nabla\phi_{\PP}(t),\PP \big)\six
=&-D_{\PP}H(\vareps(\uu_{\PP}(t)),\nabla\phi_{\PP}(t),\PP)-D_{\PP}\omega(\PP)+\lol_3(t).
\end{align*}
The Bochner-measurability of $t\mapsto\lol_3(t)$ follows directly from Assumption \ref{A4}, and $t\mapsto D_{\PP}\omega(\PP)$ is also Bochner-measurable since it is time-independent. We thus still need to verify the Bochner-measurability of the mapping $t\mapsto D_{\PP}{H}(\vareps(\uu_{\PP}(t)),\nabla\phi_{\PP}(t),\PP)$. We firstly show the Bochner-measurability of the mapping $t\mapsto (\uu_{\PP}(t),\phi_{\PP}(t))$. This means that we need to find a sequence $(\uu_n(t),\phi_n(t))$ of simple functions such that
\begin{align}\label{convergence of simple u phi}
(\uu_n(t),\phi_n(t))\to (\uu_{\PP}(t),\phi_{\PP}(t))\quad\text{ in $\W$ for a.a. $t\in[0,T]$}
\end{align}
as $n\to\infty$. Recall from \eqref{ltp} that $\lol(t,\PP)$ is equal to
\begin{align}\label{deff of l_1(PP) and l_2(t)}
\lol(t,\PP)&=\lol_1(\PP)+\lol_2(t)
\end{align}
where $\lol_1$ and $\lol_2$ are given by \eqref{ltp1} and \eqref{ltp2} respectively. Since $\lol_2$ is Bochner-measurable, there exist simple functions $\{\lol_{2}^n\}_{n\in\N}$ and some measurable set $\Theta\subset[0,T]$ with full measure such that
\begin{align}\label{converence of simple}
\lol_{2}^n(t)\to \lol_2(t)\quad\text{ in $\W^*$ for all $t\in \Theta$}
\end{align}
as $n\to\infty$. Now let $t\in\Theta$. It is immediate that ${\ell _1}(\PP)$ is an element of $\W^*$ for each $\PP\in (Y,X)_{\frac{1}{\tau},\tau}\subset \big(C(\yoverline{\Omaga})\big)^d$ by using Assumption \ref{A2}. Hence Assumption \ref{A3} is applicable with the r.h.s. function
$${\ell^n}(t,\PP):={\ell _1}(\PP)+{\ell^n _2}(t)$$
and we denote by $(\uu_n(t),\phi_{n}(t))$ the unique $\W$-solution of
\begin{align}\label{unphindef}
{\LL_{\PP}}(\uu_n(t),\phi_n(t))={\ell^n}(t,\PP).
\end{align}
Since ${\ell^n}(t,\PP)$ is a simple function, either is $(\uu_n(t),\phi_n(t))$. Moreover, we obtain that
$$ {\LL_{\PP}}(\uu_n(t)-\uu_{\PP}(t),\phi_n(t)-\phi_{\PP}(t))={\ell^n}(t,\PP)-{\ell (}t,\PP)={\ell^n_2}(t)-{\ell _2}(t).$$
This, together with \eqref{converence of simple} and \eqref{haochang4} immediately implies \eqref{convergence of simple u phi}. Next, we show that
\begin{alignat}{2}\label{convergence of simple D_pH}
D_{\PP}H(\vareps(\uu_n(t)),\nabla\phi_n(t),\PP)\to D_{\PP}H(\vareps(\uu_{\PP}(t)),\nabla\phi_{\PP}(t),\PP)& \quad&&\text{ in }X\text{ for all }t\in\Theta
\end{alignat}
as $n\to\infty$. Since $\{D_{\PP}H(\vareps(\uu_n(t)),\nabla\phi_n(t),\PP)\}_{n\in\N}$ is a sequence of simple functions, \eqref{convergence of simple D_pH} will complete the proof of Step 1. Let $\BPP\in \big(W^{1,p'}_{\Gamma_{\PP}}(\Omaga)\big)^d$. We obtain from \eqref{H} that
\begin{align}\label{expression of h}
&\int_{\Omega}D_{\PP}H(\vareps(\uu),\nabla\phi,\PP)\cdot\BPP d{\x}\nonumber\ten
=&\int_{\Omaga}\frac{1}{2}D_{\PP}\mathbb{B}(\PP)\,\BPP\begin{pmatrix}\vareps(\uu)-\vareps^0(\PP)\\ -\nabla\phi\end{pmatrix}:\begin{pmatrix}\vareps(\uu)-\vareps^0(\PP)\\ \nabla\phi\end{pmatrix}\,dx\nonumber\ten
+&\int_{\Omaga}\mathbb{B}(\PP)\,\begin{pmatrix}-D_{\PP}\vareps^0(\PP)\BPP\\ \0\end{pmatrix}:\begin{pmatrix}\vareps(\uu)-\vareps^0(\PP)\\ \nabla\phi\end{pmatrix}\,d\x+\int_{\Omaga}\nabla\phi\cdot\BPP\,d\x\nonumber\ten
=&:\,\mathcal{S}_1(\uu,\phi,\PP)+\mathcal{S}_2(\uu,\phi,\PP)+\mathcal{S}_3(\phi).
\end{align}
We firstly show that the integral in \eqref{expression of h} is well-defined for $(\uu,\phi)\in\W$. Using Assumption \ref{A2} and the fact that $\PP\in \big(C(\yoverline{\Omaga})\big)^d$ we can find some $K_{\PP}\in(1,\infty)$ such that
$$ \sup_{x\in\bar{\Omaga}}\Big\{\,\big|\mathbb{B}(\PP({\x}))\big|,\big|D_{\PP}\mathbb{B}(\PP({\x}))\big|,\big|\vareps^0(\PP({\x}))\big|,\big|D_{\PP}\vareps^0(\PP({\x}))\big|\Big\}+\|\PP\|_{L^{\infty}(\Omaga)}\leq K_{\PP}.$$
We obtain that
\begin{align}
&\left|\mathcal{S}_1(\uu,\phi,\PP)+\mathcal{S}_2(\uu,\phi,\PP)\right|\nonumber\ten
\leq&\int_{\Omaga}CK_{\PP}\left(|\vareps(\uu)|^2+|\nabla\phi|^2+K^2_{\PP}\right)|\BPP|+K^2_{\PP}\left(|\vareps(\uu)|+|\nabla\phi|+K_{\PP}\right)|\BPP|\,d\x\nonumber\ten
\leq&\int_{\Omaga}CK^3_{\PP}\left(1+|\vareps(\uu)|^2+|\nabla\phi|^2\right)|\BPP|\,d\x\nonumber\ten
\leq&\,CK^3_{\PP}\left(\|1\|^2_{L^p(\Omaga)}+\|\vareps(\uu)\|^2_{L^p(\Omaga)}+\|\nabla\phi\|^2_{L^p(\Omaga)}\right)\|\BPP\|_{L^{\frac{p}{p-2}}(\Omaga)}\nonumber\ten
\leq &\, CK^3_{\PP}\left(1+\|\uu\|^2_{W^{1,p}(\Omaga)}+\|\phi\|^2_{W^{1,p}(\Omaga)}\right)\|\BPP\|_{W^{1,p'}(\Omaga)},\label{zuihou4}
\end{align}
where we used the equality $\frac{1}{p}+\frac{1}{p}+\frac{p-2}{p}=1$ for the H\"older's inequality and the fact that $W^{1,p'}(\Omaga)\hookrightarrow L^{\frac{p}{p-2}}(\Omaga)$ for $p\in(d,\infty)$, which follows from the implication
\begin{align}\label{why p larger d}
p\in(d,\infty)\Rightarrow 1-\frac{d}{p'}\geq 0-d\big/\left(\frac{p}{p-2}\right).
\end{align}
Analogously, for $\mathcal{S}_3(\phi)$ we have
\begin{align}
\left|\mathcal{S}_3(\phi)\right|\leq C\|\nabla\phi\|_{L^p(\Omaga)}\|\BPP\|_{L^{\frac{p}{p-2}}(\Omaga)}\leq C\|\phi\|_{W^{1,p}(\Omaga)}\|\BPP\|_{W^{1,p'}(\Omaga)}.\label{zuihou5}
\end{align}
The well-definedness of the integral in \eqref{expression of h} then follows from the fact that $(\uu,\phi)\in\W\subset \big(W^{1,p}(\Omaga)\big)^d$. Now we consider the difference
\begin{align}
&\int_{\Omega}\left(D_{\PP}H(\vareps(\uu_{n}(t)),\nabla\phi_{n}(t),\PP)-D_{\PP}H(\vareps(\uu_{\PP}(t)),\nabla\phi_{\PP}(t),\PP)\right)\cdot\BPP d{\x} \nb\ten
=&\left(\SSS_1(\uu_n(t),\phi_n(t),\PP)-\SSS_1(\uu_{\PP}(t),\phi_{\PP}(t),\PP)\right)\nonumber\six
+&\left(\SSS_2(\uu_n(t),\phi_n(t),\PP)-\SSS_2(\uu_{\PP}(t),\phi_{\PP}(t),\PP)\right)\nonumber\six
+&\left(\SSS_3(\phi_n(t))-\SSS_3(\phi_{\PP}(t))\right).
\end{align}
We recall that
\begin{gather}
\LL_{\PP}(\uu_{n}(t),\phi_{n}(t))=\lol_1(\PP)+\lol^n_2(t),\label{lp est1}\\[4pt]
\LL_{\PP}(\uu_{\PP}(t),\phi_{\PP}(t))=\lol_1(\PP)+\lol_2(t),\label{lp est2}\\[4pt]
{\LL_{\PP}}(\uu_n(t)-\uu_{\PP}(t),\phi_n(t)-\phi_{\PP}(t))={\ell^n_2}(t)-{\ell _2}(t).\label{lp est3}
\end{gather}
From \eqref{ltp1} it follows
\begin{align}
\|l_1(\PP)\|_{\W^*}\leq CK^2_{\PP}.
\end{align}
Applying \eqref{haochang4} on \eqref{lp est1} to \eqref{lp est3}, we deduce that
\begin{align}
\|(\uu_{\PP}(t),\phi_{\PP}(t))\|_{\W}&\leq C C_{\PP}K^2_{\PP}(1+\|l_2(t)\|_{\W^*}),\label{un est1}\\[4pt]
\|(\uu_{n}(t),\phi_{n}(t))\|_{\W}&\leq C C_{\PP}K^2_{\PP}(1+\|l^n_2(t)\|_{\W^*}),\label{un est2}\\[4pt]
\|(\uu_{n}(t)-\uu_{\PP}(t),\phi_{n}(t)-\phi_{\PP}(t))\|_{\W}&\leq C C_{\PP}\|l_2^n(t)-l_2(t)\|_{\W^*},\label{un est3}
\end{align}
where $C_{\PP}$ is the constant given by \eqref{haochang4}. Using telescoping technique and \eqref{un est1} to \eqref{un est3} we obtain that
\begin{align}
&\left|\SSS_1(\uu_n(t),\phi_n(t),\PP)-\SSS_1(\uu_{\PP}(t),\phi_{\PP}(t),\PP)\right|\nonumber\ten
\leq & \int_{\Omaga}C K_{\PP}\left(|\vareps(\uu_n(t))-\vareps(\uu_{\PP}(t))|+|\nabla\phi_n(t)-\nabla\phi_{\PP}(t)|\right)\left|\BPP\right|\nb\\[4pt]
\cdot&\left(|\vareps(\uu_n(t))|+|\vareps(\uu_{\PP}(t))|+|\nabla\phi_n(t)|+|\nabla\phi_{\PP}(t)|+K_{\PP}\right)\,d\x\nb\ten
\leq & \int_{\Omaga}C K^2_{\PP}\left(|\vareps(\uu_n(t))-\vareps(\uu_{\PP}(t))|+|\nabla\phi_n(t)-\nabla\phi_{\PP}(t)|\right)\left|\BPP\right|\nb\\[4pt]
\cdot&\left(1+|\vareps(\uu_n(t))|+|\vareps(\uu_{\PP}(t))|+|\nabla\phi_n(t)|+|\nabla\phi_{\PP}(t)|\right)\,d\x\nb\ten
\leq &\,C K^2_{\PP}\left(\|\vareps(\uu_n(t))-\vareps(\uu_{\PP}(t))\|_{L^p(\Omaga)}+\|\nabla\phi_n(t)-\nabla\phi_{\PP}(t)\|_{L^p(\Omaga)}\right)\|\BPP\|_{L^{\frac{p}{p-2}}(\Omaga)}\nb\\[4pt]
\cdot& \left(1+\|\vareps(\uu_n(t))\|_{L^p(\Omaga)}+\|\vareps(\uu_{\PP}(t))\|_{L^p(\Omaga)}+\|\nabla\phi_n(t)\|_{L^p(\Omaga)}+\|\nabla\phi_{\PP}(t)\|_{L^p(\Omaga)}\right)\nb\ten
\leq &\,CK^4_{\PP}\,C^2_{\PP}\|l_2^n(t)-l_2(t)\|_{\W^*}\left(1+\|l_2^n(t)\|_{\W^*}+\|l_2(t)\|_{\W^*}\right)\|\BPP\|_{W^{1,p'}(\Omaga)},\label{diff est1}
\end{align}
\begin{align}
&\left|\SSS_2(\uu_n(t),\phi_n(t),\PP)-\SSS_2(\uu_{\PP}(t),\phi_{\PP}(t),\PP)\right|\nonumber\ten
\leq & \int_{\Omaga}K^2_{\PP}\left(|\vareps(\uu_n(t))-\vareps(\uu_{\PP}(t))|+|\phi_n(t)-\phi_{\PP}(t)|\right)|\BPP|\,d\x\nb\ten
\leq &\,CK^2_{\PP}\left(\|\vareps(\uu_n(t))-\vareps(\uu_{\PP}(t))\|_{L^p(\Omaga)}+\|\phi_n(t)-\phi_{\PP}(t)\|_{L^p(\Omaga)}\right)\|\BPP\|_{L^{\frac{p}{p-2}}(\Omaga)}\nb\ten
\leq &\,CK^2_{\PP}C_{\PP}\|l_2^n(t)-l_2(t)\|_{\W^*}\|\BPP\|_{W^{1,p'}(\Omaga)}\label{diff est2}
\end{align}
and
\begin{align}
&\left|\SSS_3(\phi_n(t))-\SSS_3(\phi_{\PP}(t))\right|\nb\\[4pt]
\leq&\, C\|\nabla\phi_n(t)-\nabla\phi_{\PP}(t)\|_{L^p(\Omaga)}\|\BPP\|_{L^{\frac{p}{p-2}}(\Omaga)}\nb\\[4pt]
\leq&\,CC_{\PP}\|l_2^n(t)-l_2(t)\|_{\W^*}\|\BPP\|_{W^{1,p'}(\Omaga)}.\label{diff est3}
\end{align}
From \eqref{converence of simple} we know that the norms $\|l^n_2(t)\|_{\W^*}$ are uniformly bounded for all $n\in\N$, thus \eqref{diff est1} to \eqref{diff est3} and \eqref{converence of simple} immediately imply \eqref{convergence of simple D_pH}. This completes the proof of Step 1.

\subsubsection*{Step 2: validity of condition (S)}
We recall the condition (S): Let $M>0$ be an arbitrary positive number. Then we need to show that there exists some $h_M\in L^\tau(0,T)$, such that for all $\PP_1,\PP_2\in(Y,X)_{\frac{1}{\tau},\tau}$ with
\begin{align*}
\max\left\{\,\|\PP_1\|_{(Y,X)_{\frac{1}{\tau},\tau}},\|\PP_2\|_{(Y,X)_{\frac{1}{\tau},\tau}}\,\right\}\leq M
\end{align*}
we have
\begin{align}\label{estimate stp}
\|S(t,\PP_1)-S(t,\PP_2)\|_{X}\leq h_M(t)\|\PP_1-\PP_2\|_{(Y,X)_{\frac{1}{\tau},\tau}}
\end{align}
for a.a. $t\in(0,T)$. From Lemma \ref{linfty} we deduce that the set
\begin{align}\label{def of v}
\mathcal{V}:=\{\PP\in (Y,X)_{\frac{1}{\tau},\tau}:\|\PP\|_{(Y,X)_{\frac{1}{\tau},\tau}}\leq M\}
\end{align}
is a bounded subset of $\big(L^\infty(\Omaga)\big)^d$. We may then assume that $\|\PP\|_{L^\infty(\Omaga)}\leq M'$ for all $\PP\in\mathcal{V}$ with some $M'>0$ depending on $M$. Using Assumption \ref{A2} we can find some $C_M\in(\max\{1,M'\},\infty)$ such that
\begin{multline*}
\sup_{\boldsymbol{p}\in\R^d,|\boldsymbol{p}|\leq M'}\Big\{{\mathrm{Lip}}_{D_{\boldsymbol{p}}\C},{\mathrm{Lip}}_{D_{\boldsymbol{p}}\e},{\mathrm{Lip}}_{D_{\boldsymbol{p}}\vareps^0},{\mathrm{Lip}}_{D_{\boldsymbol{p}}\eps},{\mathrm{Lip}}_{D_{\boldsymbol{p}}\omega}\\[4pt]
\left|\C(\boldsymbol{p})\right|,\left|\e(\boldsymbol{p})\right|,\left|\vareps^0(\boldsymbol{p})\right|,\left|\eps(\boldsymbol{p})\right|,\left|\omega(\boldsymbol{p})\right|\\[4pt]
\left|D_{\boldsymbol{p}}\C(\boldsymbol{p})\right|,\left|D_{\boldsymbol{p}}\e(\boldsymbol{p})\right|,\left|D_{\boldsymbol{p}}\vareps^0(\boldsymbol{p})\right|,\left|D_{\boldsymbol{p}}\eps(\boldsymbol{p})\right|,\left|D_{\boldsymbol{p}}\omega(\boldsymbol{p})\right|\Big\}\leq C_M,
\end{multline*}
where for a Lipschitz function $f$ defined on $\R^d$, ${\mathrm{Lip}}_{f}$ denotes its Lipschitz constant on the closed ball $\bar{B}_{M'}(0)\subset\R^d$. Let $\PP_1,\,\PP_2\in \mathcal{V}$. Consider the difference
\begin{align}
&S(t,\PP_1)-S(t,\PP_2)\nb\\
=&-\Big(D_{\PP}H(\vareps(\uu_{\PP_1}(t)),\nabla\phi_{\PP_1}(t),\PP_1)-D_{\PP}H(\vareps(\uu_{\PP_2}(t)),\nabla\phi_{\PP_2}(t),\PP_2)\Big)\nb\\
&-\Big(D_{\PP}\omega(\PP_1)-D_{\PP}\omega(\PP_2)\Big)\nb\\
=&:\mathcal{I}_1+\mathcal{I}_2.
\end{align}
To estimate $\mathcal{I}_2$, we obtain that
\begin{align*}
&\int_{\Omaga}|D_{\PP}\omega(\PP_1)-D_{\PP}\omega(\PP_2)|\cdot|\BPP|dx\\[4pt]
\leq &\,C\|D_{\PP}\omega(\PP_1)-D_{\PP}\omega(\PP_2)\|_{L^{p}(\Omaga)}\|\BPP\|_{L^{\frac{p}{p-2}}(\Omaga)}\\[4pt]
\leq&\, CC_M\|\PP_1-\PP_2\|_{L^{p}(\Omaga)}\|\BPP\|_{W^{1,p'}(\Omaga)}\\[4pt]
\leq& \,CC_M\|\PP_1-\PP_2\|_{L^{\infty}(\Omaga)}\|\BPP\|_{W^{1,p'}(\Omaga)}\\[4pt]
\leq &\,CC_M\|\PP_1-\PP_2\|_{(Y,X)_{\frac{1}{\tau},\tau}}\|\BPP\|_{W^{1,p'}(\Omaga)},
\end{align*}
from which we infer that
\begin{align}\label{est of I2}
\|\mathcal{I}_2\|_{X}\leq CC_M\|\PP_1-\PP_2\|_{(Y,X)_{\frac{1}{\tau},\tau}}.
\end{align}
To estimate $\mathcal{I}_1$, we use \eqref{expression of h} to obtain that
\begin{align}
&\int_{\Omega}\left(D_{\PP}H(\vareps(\uu_{\PP_1}(t)),\nabla\phi_{\PP_1}(t),\PP_1)-D_{\PP}H(\vareps(\uu_{\PP_2}(t)),\nabla\phi_{\PP_2}(t),\PP_2)\right)\cdot\BPP d{\x} \nb\ten
=&\left(\SSS_1(\uu_{\PP_1}(t),\phi_{\PP_1}(t),\PP_1)-\SSS_1(\uu_{\PP_2}(t),\phi_{\PP_2}(t),\PP_2)\right)\nonumber\six
+&\left(\SSS_2(\uu_{\PP_1}(t),\phi_{\PP_1}(t),\PP_1)-\SSS_2(\uu_{\PP_2}(t),\phi_{\PP_2}(t),\PP_2)\right)\nonumber\six
+&\left(\SSS_3(\phi_{\PP_1}(t))-\SSS_3(\phi_{\PP_2}(t))\right).
\end{align}
From the definition of $(\uu_{\PP_j},\phi_{\PP_j})$ with $j=1,2$ we know that
\begin{align}
&\LL_{\PP_j}(\uu_{\PP_i}(t),\phi_{\PP_j}(t))=\lol_1(\PP_j)+\lol_2(t)\label{pp1-pp2 1}
\end{align}
and
\begin{align}
&\LL_{\PP_1}(\uu_{\PP_1}(t)-\uu_{\PP_2}(t),\phi_{\PP_1}(t)-\phi_{\PP_2}(t))\nb\\[4pt]
=&{\lol _1}(\PP_1)-{\lol_1}(\PP_2)-\Big([\LL_{\PP_1}-\LL_{\PP_2}]\big(\uu_{\PP_2}(t),\phi_{\PP_2}(t)\big)\Big).\label{pp1-pp2 2}
\end{align}
From \eqref{ltp1} and \eqref{pp1-pp2 6} it follows
\begin{align}
\|\ell_1(\PP_j)\|_{\W^*}&\leq CC^2_M,\label{pp1-pp2 3}\\[4pt]
\|\ell_1(\PP_1)-\ell_1(\PP_2)\|_{\W^*}&\leq CC^2_M\|\PP_1-\PP_2\|_{L^\infty(\Omaga)}\label{pp1-pp2 8}
\end{align}
and
\begin{align}
&\Big\|[\LL_{\PP_1}-\LL_{\PP_2}]\big(\uu_{\PP_2}(t),\phi_{\PP_2}(t)\big)\Big\|_{\W^*}\nb\\[4pt]
\leq&\, CC_M\|\PP_1-\PP_2\|_{L^\infty(\Omaga)}\|(\uu_{\PP_2}(t),\phi_{\PP_2}(t))\|_{\W}\label{pp1-pp2 7}.
\end{align}
Applying \eqref{haochang4} on \eqref{pp1-pp2 1} and \eqref{pp1-pp2 2} and using \eqref{pp1-pp2 3} and \eqref{pp1-pp2 8} we obtain that
\begin{align}
\|(\uu_{\PP_j}(t),\phi_{\PP_j}(t))\|_{\W}\leq C C^*C^2_M(1+\|l_2(t)\|_{\W^*})\label{pp1-pp2 4}
\end{align}
with $C^*$ the number from Lemma \ref{haochang6} corresponding to the set $\mathcal{V}$ (notice here that Lemma \ref{haochang6} is applicable, since $\mathcal{V}$ is due to Lemma \ref{haochang4} a closed and bounded subset of some H\"older space, and H\"older spaces are compactly embedded into the space of continuous functions), and
\begin{align}
&\|(\uu_{\PP_1}(t)-\uu_{\PP_2}(t),\phi_{\PP_1}(t)-\phi_{\PP_2}(t))\|_{\W}\nb\six
\leq&\,CC^*C^2_M\|\PP_1-\PP_2\|_{L^\infty(\Omaga)}+CC^*C_M\|\PP_1-\PP_2\|_{L^\infty(\Omaga)}\|(\uu_{\PP_2}(t),\phi_{\PP_2}(t))\|_{\W}\nb\six
\leq&\,CC^*C^2_M\|\PP_1-\PP_2\|_{L^\infty(\Omaga)}\left(1+\|(\uu_{\PP_2}(t),\phi_{\PP_2}(t))\|_{\W}\right)\nb\six
\overset{\eqref{pp1-pp2 4}}{\leq}&\,CC^*C^2_M\|\PP_1-\PP_2\|_{L^\infty(\Omaga)}\left(C C^*C^2_M(1+\|l_2(t)\|_{\W^*})\right)\nb\six
\leq&\,C(C^*)^2C^4_M\left(1+\|l_2(t)\|_{\W^*}\right)\|\PP_1-\PP_2\|_{L^\infty(\Omaga)}.\label{pp1-pp2 5}
\end{align}
Combining with telescoping technique, \eqref{pp1-pp2 4} and \eqref{pp1-pp2 5} yield
\begin{align}
&\left|\SSS_1(\uu_{\PP_1}(t),\phi_{\PP_1}(t),\PP_1)-\SSS_1(\uu_{\PP_2}(t),\phi_{\PP_2}(t),\PP_2)\right|\nb\six
\leq& \int_{\Omaga}CC^3_M|\PP_1-\PP_2|\left(1+|\vareps(\uu_{\PP_1}(t))|^2+|\nabla\phi_{\PP_1}(t)|^2\right)|\BPP|\,d\x\nb\\[4pt]
+&\int_{\Omaga}CC^3_M\Big(|\vareps(\uu_{\PP_1}(t))-\vareps(\uu_{\PP_2}(t))|+|\nabla\phi_{\PP_1}(t)-\nabla\phi_{\PP_2}(t)|+|\PP_1-\PP_2|\Big)\nb\\[4pt]
\cdot&\,\Big(1+\sum_{j=1}^2\big(|\vareps(\uu_{\PP_j}(t))|+|\nabla\phi_{\PP_j}(t)|\big)\Big)|\BPP|\,d\x\nb\six
\leq&\,CC^3_M\|\PP_1-\PP_2\|_{L^p(\Omaga)}\left(1+\|\vareps(\uu_{\PP_1}(t))\|_{L^p(\Omaga)}^2+\|\nabla\phi_{\PP_1}(t)\|_{L^p(\Omaga)}^2\right)\|\BPP\|_{L^{\frac{p}{p-2}}(\Omaga)}\nb\\[4pt]
+&\,CC^3_M\|\BPP\|_{L^{\frac{p}{p-2}}(\Omaga)}\Big(\|\vareps(\uu_{\PP_1}(t))-\vareps(\uu_{\PP_2}(t))\|_{L^p(\Omaga)}+\|\PP_1-\PP_2\|_{L^p(\Omaga)}\nb\\[4pt]
+&\|\nabla\phi_{\PP_1}(t)-\nabla\phi_{\PP_2}(t)\|_{L^p(\Omaga)}\Big)\cdot\Big(1+\sum_{j=1}^2\big(\|\vareps(\uu_{\PP_j}(t))\|_{L^p(\Omaga)}+\|\nabla\phi_{\PP_j}(t)\|_{L^p(\Omaga)}\big)\Big)\nb\\
\leq &\,CC_M^3\left(\|(\uu_{\PP_1}(t)-\uu_{\PP_2}(t),\phi_{\PP_1}(t)-\phi_{\PP_2}(t))\|_{\W}+\|\PP_1-\PP_2\|_{L^\infty(\Omaga)}\right)\nb\\[4pt]
\cdot &\,\big(1+\sum_{j=1}^2\|\left(\uu_{\PP_j}(t),\phi_{\PP_j}(t)\right)\|_{\W}\big)\|\BPP\|_{W^{1,p'}(\Omaga)}\nb\six
\leq &\,C(C^*)^3C^9_M\left(1+\|l_2(t)\|^2_{\W^*}\right)\|\PP_1-\PP_2\|_{L^\infty(\Omaga)}\|\BPP\|_{W^{1,p'}(\Omaga)}.\label{zuihou1}
\end{align}
In the same manner, we obtain that
\begin{align}
&\left|\SSS_2(\uu_{\PP_1}(t),\phi_{\PP_1}(t),\PP_1)-\SSS_2(\uu_{\PP_2}(t),\phi_{\PP_2}(t),\PP_2)\right|\nb\six
\leq &\int_{\Omaga}CC^3_M|\PP_1-\PP_2|\left(1+|\vareps(\uu_{\PP_1}(t))|+|\nabla\phi_{\PP_1}(t)|\right)|\BPP|\,d\x\nb\six
+ &\int_{\Omaga}C^3_M\Big(|\vareps(\uu_{\PP_1}(t))-\vareps(\uu_{\PP_2}(t))|+|\nabla\phi_{\PP_1}(t)-\nabla\phi_{\PP_2}(t)|+|\PP_1-\PP_2|\Big)|\BPP|\,d\x\nb\six
\leq &\,CC^3_M\|\PP_1-\PP_2\|_{L^p(\Omaga)}\Big(1+\|\vareps(\uu_{\PP_1}(t))\|_{L^p(\Omaga)}+\|\nabla\phi_{\PP_1}(t)\|_{L^p(\Omaga)}\Big)\|\BPP\|_{L^{\frac{p}{p-2}}(\Omaga)}\nb\six
+&\,CC^3_M\Big(\|\vareps(\uu_{\PP_1}(t))-\vareps(\uu_{\PP_2}(t))\|_{L^p(\Omaga)}+\|\nabla\phi_{\PP_1}(t)-\nabla\phi_{\PP_2}(t)\|_{L^p(\Omaga)}\nb\six
+&\,\|\PP_1-\PP_2\|_{L^p(\Omaga)}\Big)\|\BPP\|_{L^{\frac{p}{p-2}}(\Omaga)}\nb\six
\leq &\,CC^*C^5_M\left(1+\|l_2(t)\|_{\W^*}\right)\|\PP_1-\PP_2\|_{L^\infty(\Omaga)}\|\BPP\|_{W^{1,p'}(\Omaga)}\nb\six
+&\,C(C^*)^2C^7_M\left(1+\|l_2(t)\|_{\W^*}\right)\|\PP_1-\PP_2\|_{L^\infty(\Omaga)}\|\BPP\|_{W^{1,p'}(\Omaga)}\nb\six
\leq &\,C(C^*)^3C^9_M\left(1+\|l_2(t)\|^2_{\W^*}\right)\|\PP_1-\PP_2\|_{L^\infty(\Omaga)}\|\BPP\|_{W^{1,p'}(\Omaga)}\label{zuihou2}
\end{align}
and
\begin{align}
&\left|\SSS_3(\phi_{\PP_1}(t))-\SSS_3(\phi_{\PP_2}(t))\right|\nb\six
\leq &\int_{\Omaga}|\nabla\phi_{\PP_1}(t)-\nabla\phi_{\PP_2}(t)|\cdot|\BPP|\,d\x\nb\six
\leq &\,C\|\nabla\phi_{\PP_1}(t)-\nabla\phi_{\PP_2}(t)\|_{L^p(\Omaga)}\|\BPP\|_{L^{\frac{p}{p-2}}(\Omaga)}\nb\six
\leq &\,C(C^*)^2C^4_M\left(1+\|l_2(t)\|_{\W^*}\right)\|\PP_1-\PP_2\|_{L^\infty(\Omaga)}\|\BPP\|_{W^{1,p'}(\Omaga)}\nb\six
\leq &\,C(C^*)^3C^9_M\left(1+\|l_2(t)\|^2_{\W^*}\right)\|\PP_1-\PP_2\|_{L^\infty(\Omaga)}\|\BPP\|_{W^{1,p'}(\Omaga)}\label{zuihou3}
\end{align}
for some sufficiently large $C$. Hence
\begin{align}\label{est of I1}
\|\mathcal{I}_1\|_{X}&\leq C(C^*)^3C^9_M\left(1+\|l_2(t)\|^2_{\W^*}\right)\|\PP_1-\PP_2\|_{L^\infty(\Omaga)}\nb\six
& \leq C(C^*)^3C^9_M\left(1+\|l_2(t)\|^2_{\W^*}\right)\|\PP_1-\PP_2\|_{(Y,X)_{\frac{1}{\tau},\tau}}.
\end{align}
From Assumption \ref{A4} it follows that $(1+\|l_2(t)\|^2_{\W^*})\in L^{\tau}(0,T)$, thus setting
$$h_M(t)=C(C^*)^3C^9_M\left(1+\|l_2(t)\|^2_{\W^*}\right)+CC_M$$
with $CC_M$ from \eqref{est of I2}, completes the proof of Step 2.
\subsubsection*{Step 3: verification of $S(\cdot,\0)\in L^\tau(0,T;X)$}
Recall that
$$S(t,\0)=-D_{\PP}H(\vareps(\uu_{\0}(t)),\nabla\phi_{\0}(t),\0)-D_{\PP}\omega(\0)+\lol_3(t).$$
From Assumption \ref{A4} we know that $\lol_3\in L^\tau(0,T;X)$; Since $-D_{\PP}\omega(\0)$ is a constant, it is also an element of $L^\tau(0,T;X)$; Finally, by setting $(\uu,\phi)=(\uu_{\0}(t),\phi_{\0}(t))$ in \eqref{zuihou4} and \eqref{zuihou5} and using \eqref{un est1} with $\PP=\0$ and Assumption \ref{A4} we conclude that
$$-D_{\PP}H(\vareps(\uu_{\0}(t)),\nabla\phi_{\0}(t),\0)\in L^\tau(0,T;X).$$
Hence we infer that $S(\cdot,\0)\in L^\tau(0,T;X)$.

Sum up all, we see that all conditions of Theorem \ref{pruss2002} are satisfied and one obtains from Theorem \ref{pruss2002} a unique local in time solution
\begin{align*}
\PP\in W^{1,r}(0,\hat{T};\big(W^{-1,p}_{\Gamma_{\PP}}(\Omaga)\big)^d)\cap L^{r}(0,\hat{T};\big(W^{1,p}_{\Gamma_{\PP}}(\Omaga)\big)^d)
\end{align*}
of the equation
\begin{align*}
\PP'(t)-\Delta\PP(t)=S\big(t,\PP(t)\big)
\end{align*}
on the maximal time interval $(0,\hat{T})$ for some $\hat{T}\in(0,T]$. For the given pair $(t,\PP(t))$ one obtains the unique weak solution $(\uu(t),\phi(t))$ of \eqref{elliptic part}, that is
\begin{align}
L_{\PP(t)}(\uu(t),\phi(t))={\ell (}t,\PP(t)).
\end{align}
From Proposition \ref{reg of polar} we find some positive constant $M''$ and some $\delta\in(0,1)$ such that
\begin{align}\label{M prime prime}
\max_{t\in[0,T]}\|\PP(t)\|_{C^\delta(\bar{\Omaga})}\leq M''.
\end{align}
We then define the set $\mathcal{M}$ in Lemma \ref{haochang6} by
$$\mathcal{M}:=\{\PP\in \big(C^\delta(\yoverline{\Omaga})\big)^d:\|\PP\|_{C^\delta(\bar{\Omaga})}\leq M''\} $$
and denote by $(C^*)'$ the corresponding uniform bound of inverse norm $\mathrm{L}^{-1}_{\PP}$ given in Lemma \ref{haochang6}. Using the definition of ${\ell _1}(\PP)$ given in \eqref{ltp1}, together with Assumption \ref{A2} and \eqref{M prime prime}, it follows immediately that
$$ \max_{t\in[0,T]}\|{\ell _1}(\PP(t))\|_{\W^*}\leq C_{M''}$$
for some positive constant $C_{M''}$ which only depends on $M''$. Using \eqref{haochang4} we infer that
\begin{align}\label{bound of u phi}
\|(\uu(t),\phi(t))\|_{\W}\leq (C^*)'(C_{M''}+\|{\ell _2}(t)\|_{\W^*}).
\end{align}
On the other hand, the Bochner-measurability of $(\uu,\phi):[0,T]\to\W$ can be similarly shown as by Step 1 given previously, we omit the details here. Hence the claimed regularity of $(\uu,\phi)$ given in Definition \ref{solution concept} follows immediately from Assumption \ref{A4} and \eqref{bound of u phi}. This completes the desired proof of Theorem \ref{C^1 domain}.
\end{proof}
\section{Proof of the elliptic regularity results}\label{proof of thm 2}
\subsection{Proof of Theorem \ref{2d mixed}}
In this section we give the proof of Theorem \ref{2d mixed}. The main ingredient is Theorem 6.2 from \cite{HDJonssonKneesRehberg2016}, which will be restated in Theorem \ref{higherorder} below. To introduce the result we firstly fix some notation. Let $m\in\N$. For $1\leq i\leq m$ and $\gggg_i\subset\partial\Omaga$ we define the space
\begin{equation*}
\mathbb{W}_{\Gamma}^{1,p}:=\prod_{i=1}^m W_{\gggg_i}^{1,p}(\Omaga).
\end{equation*}
The dual space of $\mathbb{W}_{\Gamma}^{1,p'}$ is denoted by $\mathbb{W}_{\Gamma}^{-1,p}$. For $\mathbb{A}\in L^\infty\big(\Omaga,\mathrm{Lin}(\R^{m}\times \R^{md},\R^{m}\times \R^{md})\big)$ we define the operator $\mathcal{A}$ by
\begin{align*}
\langle\mathcal{A}(\vv_1),\vv_2\rangle:=\int_{\Omaga}\mathbb{A}\begin{pmatrix}\vv_1\\ \nabla\vv_1\end{pmatrix}:\begin{pmatrix}\vv_2\\ \nabla\vv_2\end{pmatrix}\,d\x
\end{align*}
for $\vv_1,\vv_2\in \mathbb{W}_{\Gamma}^{1,2}$. Additionally we will also need the following assumption:
\begin{assumption}\label{assofextension}
Let $\Gamma\subset\partial\Omaga$ be closed. There exists a linear, continuous extension operator $\mathbb{E} : W^{1,1}_{{\Gamma}}(\Omaga)\rightarrow W_{\Gamma}^{1,1}(\R^d)$ which simultaneously defines a continuous extension operator $\mathbb{E} : W^{1,p}_{{\Gamma}}(\Omaga)\rightarrow W_{\Gamma}^{1,p}(\R^d)$ for all $p\in(1,\infty)$.
\end{assumption}
Having all the preliminaries we are now ready to state Theorem \ref{higherorder}:
\begin{theorem}[{\protect\cite[Thm. 6.2]{HDJonssonKneesRehberg2016}}]\label{higherorder}
Let $\Omega\subset\R^d$, $d\geq 2$, be a bounded domain. Let all $\gggg_i\subset\partial\Omega$ be $(d-1)$-sets and let the Assumption \ref{assofextension} be satisfied for all ${\Gamma}_i$. Furthermore, assume that there exists a positive constant $\eta$ such that
\begin{align}\label{coercive condition}
\langle\mathcal{A}(\vv),\vv\rangle\geq \eta\|\vv\|_{H^1(\Omaga)}
\end{align}
for all $\vv\in\mathbb{W}_{\Gamma}^{1,2}$. Then there exists some $p^*>2$ such that for all $p\in[2,p^*]$, the operator $\mathrm{L}_{\mathcal{A}}:\vv\mapsto\langle\mathcal{A}(\vv),\cdot\rangle$ defines a topological isomorphism from $\mathbb{W}_{{\Gamma}}^{1,p}$ to $\mathbb{W}_{{\Gamma}}^{-1,p}$. In particular, the norm $\|\mathrm{L}_{\mathcal{A}}^{-1}\|_{L(\mathbb{W}_{{\Gamma}}^{-1,p},\mathbb{W}_{{\Gamma}}^{1,p})}$ depends only on $\eta$ and $\|\mathbb{A}\|_{L^\infty(\Omaga)}$ and is uniform for all $p\in[2,p^*]$.
\end{theorem}
We now give the proof of Theorem \ref{2d mixed}:
\begin{proof}[Proof of Theorem \ref{2d mixed}]
To shed light on the connection, we set the space $\mathbb{W}_{\Gamma}^{1,p}$ by
$$\mathbb{W}_{\Gamma}^{1,p}=\W $$
and the operator $\mathcal{A}$ and tensor $\mathbb{A}$ by the following implicit form
\begin{align}
&\langle\mathcal{A}((\uu_1,\phi_1)),(\uu_2,\phi_2)\rangle\nb\\[5pt]
=&\int_{\Omaga}\mathbb{A}\begin{pmatrix}\uu_1\\\phi_1\\ \nabla\uu_1\\ \nabla\phi_1\end{pmatrix}:\begin{pmatrix}\uu_2\\\phi_2\\ \nabla\uu_2\\ \nabla\phi_2\end{pmatrix}\,d\x \nb\\[5pt]
=&\,\LL_{\PP}(\uu_1,\phi_1)[(\uu_2,\phi_2)].
\end{align}
It suffices to show that the conditions in Theorem \ref{higherorder} are satisfied. Firstly we point out that the Assumption \ref{B1} already implies Assumption \ref{assofextension}, see \cite[Thm. 4.15]{HDJonssonKneesRehberg2016}. On the other hand, from Assumption \ref{B2} we know that the tensor $\mathbb{A}$ is of class $L^\infty$. It is therefore only left to show the coercivity condition \eqref{coercive condition}. From Assumption \ref{B3}, Poincar\'e's and Korn's inequalities and the fact that the antidiagonals of the tensor $\mathbb{B}(\PP)$ given by \eqref{not symmetric tentor} are composed by $\e(\PP)^T$ and $-\e(\PP)$ we obtain that
\begin{align*}
&\LL_{\PP}(\uu,\phi)[(\uu,\phi)]\nb\\
=&\int_{\Omaga}\C(\PP)\vareps(\uu):\vareps(\uu)+\eps(\PP)\nabla\phi\cdot\nabla\phi\,d\x\nb\\
\geq &\,\alpha\left(\|\vareps(\uu)\|^2_{L^2(\Omaga)}+\|\nabla\phi\|^2_{L^2(\Omaga)}\right)\nb\\
\geq &\,\alpha \min\{c_{P},c_{K}\}\left(\|\uu\|^2_{H^1(\Omaga)}+\|\phi\|^2_{H^1(\Omaga)}\right)
\end{align*}
with $\alpha>0$ from Assumption \ref{B3} and $c_{P},c_{K}>0$ the corresponding Poincar\'e's and Korn's constants, and \eqref{coercive condition} is fulfilled. This completes the desired proof.
\end{proof}
\subsection{Proof of Theorem \ref{cube or C^1 domain}}\label{proof of thm cube or c1 domain}
In this section we prove Theorem \ref{cube or C^1 domain}. In fact, we will prove the following stronger version of Theorem \ref{cube or C^1 domain}, in the sense that $p$ can be arbitrarily chosen from the interval $[2,\infty)$ and $\PP$ is merely assumed to be an element of $\big(C(\yoverline{\Omaga})\big)^d$.
\begin{theorem}\label{piezoreg}
Let the Assumptions \ref{C1} and \ref{C2} be satisfied. Let $\PP\in \big(C(\yoverline{\Omaga})\big)^d$ and $p\in[2,\infty)$. Then the differential operator $\LL_{\PP}$ defines a topological isomorphism from $\big(W_0^{1,p}(\Omaga)\big)^d$ to $\big(W^{-1,p}(\Omaga)\big)^d$.
\end{theorem}
\begin{proof}
Since $\PP$ is continuous on $\yoverline{\Omaga}$ and $\alpha$ is continuous and non-negative on $\R^d$ due to Assumption \ref{C2}, we infer that $\alpha\big(\PP(\cdot)\big)$ is continuous and non-negative on the compact set $\yoverline{\Omaga}$, thus also bounded below by some positive constant on $\yoverline{\Omaga}$. We can hence define
$$\alpha_{\PP}:=\min\left\{\alpha\big(\PP(x)\big):x\in\yoverline{\Omaga}\,\right\}>0.$$
From Assumption \ref{C2} and Poincar\'e's and Korn's inequalities we obtain that
\begin{align*}
&\LL_{\PP}(\uu,\phi)[(\uu,\phi)]\nb\\
=&\int_{\Omaga}\C(\PP)\vareps(\uu):\vareps(\uu)+\eps(\PP)\nabla\phi\cdot\nabla\phi\,d\x\nb\\
\geq&\int_{\Omaga}\alpha(\PP)(|\vareps(\uu)|^2+|\nabla\phi|^2)\,d\x\nb\\
\geq &\,\alpha_{\PP}\left(\|\vareps(\uu)\|^2_{L^2(\Omaga)}+\|\nabla\phi\|^2_{L^2(\Omaga)}\right)\nb\\
\geq &\,\alpha_{\PP} \min\{c_{P},c_{K}\}\left(\|\uu\|^2_{H^1(\Omaga)}+\|\phi\|^2_{H^1(\Omaga)}\right)
\end{align*}
and therefore ${\LL_{\PP}}$ is coercive. Since $\mathbb{B}$ is a continuous function on $\R^d$ due to Assumption \ref{C2} and $\PP\in \big(C(\yoverline{\Omaga})\big)^d$, $\mathbb{B}\big(\PP(x)\big)$ is uniformly bounded for all $x\in\yoverline{\Omaga}$ and the boundedness of ${\LL_{\PP}}$ follows from the H\"older's inequality. Thus the existence and uniqueness of a $H_0^1$-weak solution $(\uu,\phi)$ follow directly from Lax-Milgram.

It remains to show the $L^p$-isomorphism property of $\LL_{\PP}$ for arbitrary $p\in[2,\infty)$. Since $\PP(\cdot)$ is uniformly continuous on $\yoverline{\Omaga}$ and the coefficient tensor $\mathbb{B}(\cdot)$ is continuous in $\PP\in\R^d$ due to Assumption \ref{C2}, the coefficient tensor $\mathbb{B}\big(\PP(\cdot)\big)$ will also be uniformly continuous on $\yoverline{\Omaga}$. On the other hand, the coerciveness of ${\LL_{\PP}}$ implies that ${\LL_{\PP}}$ is strongly elliptic in the sense of \cite[Def. 3.36]{GiaquintaMart2012}, see for instance \cite[Thm. 4.6]{mclean2000} for a proof. Then the claim follows from the elliptic regularity results given in \cite{Acq1992} (when $\Omaga$ is a 2D-rectangle or a 3D-cuboid) and \cite{DolzmannMueller1995} (when $\Omega$ has $C^1$-boundary) for strongly elliptic systems with uniformly continuous coefficients.
\end{proof}

\subsection{Proof of Theorem \ref{poly thm}}
In this section we give the proof of Theorem \ref{poly thm}. In fact, we will show a slightly stronger version of Theorem \ref{poly thm}, in the sense that we only require that $\PP$ vanishes on $\mathrm{Sing}_{\Omaga}\subset\Gamma_{\PP}$.
\begin{theorem}\label{Lamelaplacereg}
Let the Assumption \ref{C2}, Assumption \ref{sing} and \eqref{lame and laplace} be satisfied. Let also $\Gamma_{\uu}=\Gamma_{\phi}=\partial\Omaga$. Then there exists some $p\in(3,6]$ such that ${\LL_{\PP}}$ is a topological isomorphism from $\big(W^{1,\mathfrak{p}}_0(\Omaga)\big)^4$ to $\big(W^{-1,\mathfrak{p}}(\Omaga)\big)^4$ for all $\mathfrak{p}\in[2,p]$ and for all $\PP\in \big(C(\yoverline{\Omaga})\big)^3$ with $\PP|_{\mathrm{Sing}_{\Omaga}}=\0$, where the set $\mathrm{Sing}_{\Omaga}$ is defined by Assumption \ref{sing}.
\end{theorem}
The rest of this section is devoted to the proof of Theorem \ref{Lamelaplacereg}. We firstly introduce the elliptic regularity result from \cite{MazyaVladimirRossman2010}, stated in Theorem \ref{Lamelaplacereg1} below. Based on this result, we are then able to exploit perturbation arguments to complete the proof of Theorem \ref{Lamelaplacereg}.
\begin{theorem}[\protect{\cite[Thm. 4.3.5]{MazyaVladimirRossman2010}}]\label{Lamelaplacereg1}
Let $\Omaga$ be a polyhedral domain and $\lambda,\mu,\gamma$ be given positive constants. Assume also that for each vertex ${\x}$ the corresponding cone $\mathcal{K}_{{\x}}$ is a Lipschitz graph with edges, i.e. $\mathcal{K}_{{\x}}$ has the representation ${\x}_3>\psi({\x}_1,{\x}_2)$ in Cartesian coordinates with some Lipschitz function $\psi$. Define the differential operator ${\mathrm{L}}$ by
\begin{align*}
&{\LL}:\big(H_0^1(\Omaga)\big)^4\to\big(H^{-1}(\Omaga)\big)^4,\\
&{\LL}(\uu,\phi)[\buu,\bphi]:=\int_{\Omaga}\lambda \mathrm{tr}\big(\vareps(\uu)\big)\mathrm{tr}\big(\vareps(\buu)\big)+2\mu\vareps(\uu):\vareps(\buu)+\gamma\nabla\phi\cdot\nabla\bphi d{\x}
\end{align*}
for $(\uu,\phi),\,(\buu,\bphi)\in \big(H_0^1(\Omaga)\big)^4 $. Then there exists some $\tilde{p}\in(3,\infty)$ such that ${\LL}$ is a topological isomorphism from $\big(W_0^{1,\mathfrak{p}}(\Omaga)\big)^4$ to $\big(W^{-1,\mathfrak{p}}(\Omaga)\big)^4$ for all $\mathfrak{p}\in[2,\tilde{p}]$.
\end{theorem}
Having all the preliminaries we are in the position to prove Theorem \ref{Lamelaplacereg}:
\begin{proof}[Proof of Theorem \ref{Lamelaplacereg}]
We follow a perturbation argument given by \cite[Lem. A.18]{Krumbiegel2013} to prove the statement. Consider a point ${\x}\in\yoverline{\Omaga}$. Let $\LL^0_{\PP}$ be the differential operator given by
\begin{align*}
\LL^0_{\PP}(\uu,\phi):=-\mathrm{Div}\left(\mathbb{B}\big(\PP({\x}_0)\big)\begin{pmatrix}\vareps\big(\uu\big) \\ \nabla\phi\end{pmatrix}\right)
\end{align*}
with fixed constant coefficient tensor $\mathbb{B}\big(\PP({\x}_0)\big)$ at ${\x}_0={\x}$. Firstly we show that for each ${\x}\in\yoverline{\Omaga}$, there exist some neighborhood ${U}_{{\x}}$ of ${\x}$ and some $p_{{\x}}\in(3,6]$ such that for all $\mathfrak{p}\in[2,p_{{\x}}]$ we have
\begin{align}\label{claim2}
&\LL^0_{\PP}\text{ is a topological isomorphism from}\nonumber\\
&\big(W_0^{1,\mathfrak{p}}({U}_{{\x}}\cap\Omaga)\big)^4 \text{ to } \big(W^{-1,\mathfrak{p}}({U}_{{\x}}\cap\Omaga)\big)^4.
\end{align}
We discuss different cases:
\begin{itemize}
\item If ${\x}$ is in $\Omaga$, then \eqref{claim2} follows from Theorem \ref{piezoreg}, since ${U}_{{\x}}\subset\subset\Omaga$ can be chosen as an open ball.

\item If ${\x}$ is on an open face, then by Definition \ref{polyhedral domain}, ${\x}$ is on a smooth open two-dimensional manifold. Thus we know that there exist some neighborhood ${U}_{{\x}}$ of ${\x}$ and a $C^\infty$-diffeomorphism ${\iota}$ such that ${\iota}$ maps ${U}_{{\x}}\cap\Omaga$ onto an open cube $\mathcal{C}$. According to the classical transformation arguments (see for instance \cite[Thm. 4.14]{GiaquintaMart2012}), \eqref{claim2} is equivalent to the $W_0^{1,\mathfrak{p}}\leftrightarrow W^{-1,\mathfrak{p}}$ isomorphism property of a differential operator with suitable transformed coefficient tensor $\tilde{\mathbb{B}}$, which is still strongly elliptic, on $\mathcal{C}$. But the latter claim follows immediately from Theorem \ref{piezoreg}.

\item If ${\x}\in\mathrm{Sing}^2_{\Omaga}\cup\mathrm{Sing}^3_{\Omaga}$, where $\mathrm{Sing}^2_{\Omaga},\mathrm{Sing}^3_{\Omaga}$ are defined in Assumption \ref{sing}, then by the definition of $\mathrm{Sing}^2_{\Omaga}$ and $\mathrm{Sing}^3_{\Omaga}$ we know that there exists a neighborhood ${U}_{{\x}}$ such that the intersection ${U}_{{\x}}\cap\Omaga$ is $C^\infty$-diffeomorph to an open cube. Consequently, the claim follows from Theorem \ref{piezoreg}, by using the similar transformation arguments as the ones for the points on a face given previously.

\item If ${\x}\in\mathrm{Sing}_{\Omaga}$, where $\mathrm{Sing}_{\Omaga}$ is defined in Assumption \ref{sing}, then due to \eqref{lame and laplace}, the operator $\LL^0_{\PP}$ reduces to the operator $\LL$ defined in Theorem \ref{Lamelaplacereg1}. In this case we simply take an arbitrary $U_x$ with $\yoverline{\Omaga}\subset\subset U_x$ and the claim follows from Theorem \ref{Lamelaplacereg1}.
\end{itemize}
This completes the proof of \eqref{claim2}. Now let $\PP\in \big(C(\yoverline{\Omaga})\big)^3,\ell\in \big(W^{-1,\mathfrak{p}}(\Omaga)\big)^4$ with $\mathfrak{p}\in[2,6]$ and let $(\uu,\phi)$ be the unique $H_0^1$-solution of
\begin{alignat*}{2}
&{\LL}_{\PP}(\uu,\phi)={\ell}&&\quad\text{in }\Omaga,
\end{alignat*}
whose existence and uniqueness are guaranteed by Lax-Milgram. Consider a point ${\x}\in\yoverline{\Omaga}$. Let $\eta:\R^3\to\R$ be a real valued smooth function such that $\mathrm{supp}(\eta)\subset{U}_{{\x}}$, where ${U}_{{\x}}$ is as defined in \eqref{claim2}. We also assume that $\mathfrak{p}\in[2, p_{{\x}}]$ with $p_{{\x}}$ given by \eqref{claim2}. Define $\ww:=(\uu,\phi)$$, $${V}:=(\buu,\bphi)$ and rewrite ${\LL}_{\PP}(\uu,\phi)[\buu,\bphi]$ as\footnote{We use here and below the Einstein summation notation.}
\begin{align}\label{einstein1}
{\LL}_{\PP}(\uu,\phi)[\buu,\bphi]=\int_{\Omaga}{\mathbb{A}}_{ij}^{\alpha\beta}\partial_\beta\ww^{j}\partial_\alpha{V}^{i}d{\x}
\end{align}
with coefficient tensor ${\mathbb{A}}$. For ${V}\in \big(C_0^\infty({U}_{{\x}}\cap\Omaga)\big)^4$ we obtain that
\begin{align}\label{p4}
&\int_{{U}_{{\x}}\cap\Omaga}{\mathbb{A}}_{ij}^{\alpha\beta}\partial_\beta(\eta\ww)^{j}\partial_\alpha{V}^{i}\,d{\x}\nonumber\\
=&\int_{{U}_{{\x}}\cap\Omaga}{\mathbb{A}}_{ij}^{\alpha\beta}(\partial_\beta\eta) \ww^{j}\partial_\alpha{V}^{i}+{\mathbb{A}}_{ij}^{\alpha\beta}\eta\partial_\beta\ww^{j}\partial_\alpha{V}^{i}\,d{\x}\nonumber\\
=&\int_{{U}_{{\x}}\cap\Omaga}{{\mathbb{A}}_{ij}^{\alpha\beta}(\partial_\beta\eta) \ww^{j}\partial_\alpha{V}^{i}}+{{\mathbb{A}}_{ij}^{\alpha\beta}\partial_\beta\ww^{j}\partial_\alpha(\eta{V})^{i}}+{\Big(-{\mathbb{A}}_{ij}^{\alpha\beta}\partial_\beta\ww^{j}(\partial_\alpha\eta){V}^{i}\Big)}\,d{\x}\nonumber\\
=&:\fff_a({V})+\fff(\eta{V})+\fff_b({V})=:\tilde{\fff}({V}),
\end{align}
where $\fff(\eta{V})$ is understood as the evaluation of $\fff$ at the extension of $\eta{V}$ on whole $\Omaga$ such that the extension is zero on $\Omaga\setminus{U}_{{\x}}$. From fundamental calculus one obtains that
\begin{align*}
&q\in [2,6]\Rightarrow 1-\frac{3}{2}\geq 0-\frac{3}{q},\\
&q\in[2,6]\Rightarrow 0-\frac{3}{2}\geq -1-\frac{3}{q}.
\end{align*}
Consequently, using Sobolev embedding one verifies that for $q\in[2,6]$, $\mathrm{Div}(\uu,\phi)$ and $\nabla(\uu,\phi)$ are of class $W^{-1,q}$ if $(\uu,\phi)$ is of class $H^1$. Thus $\fff_a,\fff_b$ are of class $W^{-1,q}$ for $q\in[2,6]$. Moreover, we infer that $\tilde{\fff}$ is of class $W^{-1,q}$ if and only if $\fff$ is of class $W^{-1,q}$ for $q\in[2,6]$. Since $\PP$ is uniformly continuous on $\yoverline{\Omaga}$, the coefficient tensor ${\mathbb{A}}$ is also uniformly continuous and consequently uniformly bounded on $\yoverline{\Omaga}$. Letting $q=\mathfrak{p}$ it follows from H\"older's inequality that
\begin{align}\label{regu inequ}
&\|\tilde{\fff}\|_{W^{-1,\mathfrak{p}}({U}_{{\x}}\cap\Omaga)}\nonumber\\
\leq& \,C\left(\|\fff\|_{W^{-1,\mathfrak{p}}(\Omaga)}+\|(\uu,\phi)\|_{L^{\mathfrak{p}}(\Omaga)}+\|\nabla(\uu,\phi)\|_{L^2(\Omaga)}\right).
\end{align}
Additionally we obtain that $\eta(\uu,\phi)$ is the $H_0^1$-solution of
\begin{alignat*}{2}
&{\LL}_{\PP}\big(\eta(\uu,\phi)\big)=\tilde{\fff}&&\quad\text{in }{U}_{{\x}}\cap\Omaga.
\end{alignat*}
For $r_x>0$ with $\bar{B}_{r_{{\x}}}({\x})\subset {U}_{{\x}}$ define the tensor $\mathbb{B}(\PP)^*$ on ${U}_{{\x}}\cap\Omaga$ by
\begin{equation*}
\mathbb{B}(\PP)^*(\yy):=\left\{
             \begin{array}{lr}
             \mathbb{B}\big(\PP({\x})\big), &\yy\in({U}_{{\x}}\cap\Omaga)\setminus B_{r_{{\x}}}({\x}),  \six
             \mathbb{B}\big(\PP(\yy)\big), &\yy\in({U}_{{\x}}\cap\Omaga)\cap B_{r_{{\x}}}({\x}).
             \end{array}
\right.
\end{equation*}
Denote by $\LL^*_{\PP}$ the differential operator given by
\begin{align*}
\LL^*_{\PP}(\uu,\phi):=-\mathrm{Div}\left(\mathbb{B}(\PP)^*\begin{pmatrix}\vareps\big(\uu\big) \\ \nabla\phi\end{pmatrix}\right).
\end{align*}
Using the fact that the coefficient tensor $\mathbb{B}$ is continuous on $\R^3$ from Assumption \ref{C2} and that $\PP$ is uniformly continuous on the whole $\yoverline{\Omaga}$ we deduce that $\mathbb{B}(\PP(\cdot))$ is uniformly continuous on $\yoverline{\Omaga}$. This, together with H\"older's inequality, implies that there exists a sufficiently small $r_{{\x}}>0$ such that
\begin{align}
&\|\LL^*_{\PP}-\LL^0_{\PP}\|_{L(W_0^{1,\mathfrak{p}}({U}_{{\x}}\cap\Omaga),W^{-1,\mathfrak{p}}({U}_{{\x}}\cap\Omaga))}\|(\LL^0_{\PP})^{-1}\|_{L(W^{-1,\mathfrak{p}}({U}_{{\x}}\cap\Omaga),W_0^{1,\mathfrak{p}}({U}_{{\x}}\cap\Omaga))}\nonumber\\
\leq& \,C\|\mathbb{B}\big(\PP(\cdot)\big)-\mathbb{B}\big(\PP({\x})\big)\|_{L^\infty(B_{r_{{\x}}}({\x})\cap\Omaga)}\|(\LL^0_{\PP})^{-1}\|_{L(W^{-1,\mathfrak{p}}({U}_{{\x}}\cap\Omaga),W_0^{1,\mathfrak{p}}({U}_{{\x}}\cap\Omaga))}\nonumber\\
<&\,1.
\end{align}
Combining with the small perturbation theorem \cite[Chap. 4, Thm. 1.16]{Kato1995} and \eqref{claim2} we infer that $\LL^*_{\PP}$ is a topological isomorphism from $\big(W_0^{1,\mathfrak{p}}({U}_{{\x}}\cap\Omaga)\big)^4$ to $\big(W^{-1,\mathfrak{p}}({U}_{{\x}}\cap\Omaga)\big)^4$ for all $\mathfrak{p}\in[2,p_{{\x}}]$. If we additionally let $\mathrm{supp}(\eta)\subset B_{r_{{\x}}({\x})}$, then we deduce that $\eta(\uu,\phi)$ is the $H_0^1$-solution of
\begin{alignat*}{2}
&\LL^*_{\PP}\big(\eta(\uu,\phi)\big)=\tilde{\fff}&&\quad\text{in }{U}_{{\x}}\cap\Omaga.
\end{alignat*}
This shows that $\eta(\uu,\phi)$ is of class $W^{1,\mathfrak{p}}$ on ${U}_{{\x}}\cap\Omaga$ for all $\mathfrak{p}\in[2,p_x]$.

Finally, we use a covering argument to finish the desired proof. Let ${U}_{{\x}},p_{{\x}}$ be as defined in \eqref{claim2} and $r_x$ the small number constructed in the previous step. Since $\Omaga$ is a bounded domain in $\R^3$, we can find some $k\in\N$ such that $\yoverline{\Omaga}\subset \cup_{i=1}^k B_{r_{{\x}_i}}({\x}_i)$ with ${\x}_i\in\yoverline{\Omaga}$. Define
\begin{align*}
p:=\min\left\{p_{{\x}_1},...,p_{{\x}_k}\right\}\in(3,6].
\end{align*}
At this point we emphasize that $p$ can be chosen in a (possibly non-optimal) way that $p$ is independent of the choice of $\PP$. Indeed, along the proof of \eqref{claim2} it is easy to see that one can always set $p_x\equiv\min\{\tilde{p},6\}$, where $\tilde{p}$ is given by Theorem \ref{Lamelaplacereg1} which is independent of $\PP$. Now let $\mathfrak{p}\in[2,p]$. Let $\{\eta_i\}_{i=1}^k$ be a partition of unity subordinated to the balls $\{B_{r_{{\x}_i}}({\x}_i)\}_{i=1}^k$. Then there is a corresponding $\tilde{\fff}_{i}$ of class $W^{-1,\mathfrak{p}}$, similarly defined as in \eqref{p4}, such that $\eta_i(\uu,\phi)$ is the $H_0^1$-solution of
\begin{alignat*}{2}
&\LL^*_{\PP}\big(\eta_i(\uu,\phi)\big)=\tilde{\fff}_{i}&&\quad\text{in }{U}_{{\x}_i}\cap\Omaga.
\end{alignat*}
We obtain that
\begin{align}
&\|(\uu,\phi)\|_{W_0^{1,\mathfrak{p}}(\Omaga)}\nonumber\\
\leq&\,\sum_{i=1}^k\|\eta_i(\uu,\phi)\|_{W_0^{1,\mathfrak{p}}({U}_{{\x}_i}\cap\Omaga)}\nonumber\\
\leq&\,C\sum_{i=1}^k\|\tilde{\fff}_{i}\|_{W^{-1,\mathfrak{p}}({U}_{{\x}_i}\cap\Omaga)}\nonumber\\
\leq&\,C\left(\|\fff\|_{W^{-1,\mathfrak{p}}(\Omaga)}+\|(\uu,\phi)\|_{L^{\mathfrak{p}}(\Omaga)}+\|\nabla(\uu,\phi)\|_{L^2(\Omaga)}\right)<\infty.
\end{align}
This completes the proof of Theorem \ref{Lamelaplacereg}.
\end{proof}

\subsubsection*{Acknowledgments}
This research was supported by the doctoral program of the University of Kassel. The author wants to extend his sincere thanks to Prof. Dr. Dorothee Knees for her meticulous care and great support during the author's doctoral study at the University of Kassel. The author also wants to thank the anonymous referees for their valuable comments and suggestions. Special thanks also to a referee for suggesting the replacement of the volume-preserving condition in an older version of the manuscript to the more general Ahlfors-David-regular condition, which greatly improves the applicability of the results in the present paper.

\end{document}